\documentclass[11pt,reqno]{amsart}
\usepackage{mathrsfs}
\usepackage{cite}
\usepackage{stmaryrd}
\usepackage{amssymb}
\usepackage{amsfonts}

\newcommand{\rmv}[1]{}

\def\pv#1{\ensuremath{{\bf#1}}}
\def\ZZ{\ensuremath{\mathbb Z}}

\def\inv{^{-1}}
\def\p{\varphi}

\def\J{\mathrel{{\mathscr J}}} 
\def\D{\mathrel{{\mathscr D}}} 

\def\R{\mathrel{{\mathscr R}}} 
\def\L{\mathrel{{\mathscr L}}} 

\def\e<{\leq _{E}}

\def\ov#1{\ensuremath{\overline {#1}}}

\def\malce{\mathop{\hbox{$\bigcirc$\kern-8.5pt\raise1pt
\hbox{\scriptsize$m$}\kern1.5pt}}}

\def\dom#1{\mathrm {dom}(#1)}
\def\ran#1{\mathrm {ran}(#1)}
\def\rk#1{\mathrm {rk}(#1)}

\def\1sk{^{(1)}}

\def\to{\rightarrow}

\def\pow#1{[\![#1]\!]}
\def\sgn#1#2{\mathrm{sgn}_{#2}(#1)}

\def\data{\ifcase\month\or January\or February \or March\or April\or May
\or June\or July\or August\or September\or October\or November \or
December\fi\space\number\day, \number\year}

\def\Thmname{Theorem}
\def\Propname{Proposition}
\def\Lemmaname{Lemma}
\def\Definitionname{Definition}
\newtheorem{Thm}{\Thmname}[section]
\newtheorem{Prop}[Thm]{\Propname}
\newtheorem{Lemma}[Thm]{\Lemmaname}

\newtheorem{Cor}[Thm]{Corollary}

\numberwithin{equation}{section}

\title[M\"obius Functions
  and Semigroup Representation Theory II]{M\"obius Functions
  and Semigroup Representation Theory II: \\ Character formulas and multiplicities}
\author{Benjamin Steinberg}
\address{School of Mathematics and Statistics, Carleton
University, 1125 Colonel By Drive, Ottawa, Ontario K1S 5B6,
Canada}
\thanks{The author gratefully acknowledges the support of NSERC} \email{bsteinbg@math.carleton.ca}

\date{April 20, 2006, Revised November 21, 2007}



\begin{document}

\begin{abstract}
We generalize the character formulas for multiplicities of
irreducible constituents from group theory to semigroup theory using
Rota's theory of M\"obius inversion.  The technique works for a
large class of semigroups including: inverse semigroups, semigroups
with commuting idempotents, idempotent semigroups and semigroups
with basic algebras.  Using  these tools we are able to give a
complete description of the spectra of random walks on finite
semigroups admitting a faithful representation by upper triangular
matrices over the complex numbers. These include the random walks on
chambers of hyperplane arrangements studied by Bidigare, Hanlon,
Rockmere, Brown and Diaconis.  Applications are also given to
decomposing tensor powers and exterior products of rook matrix
representations of inverse semigroups, generalizing and simplifying
earlier results of Solomon for the rook monoid.
\end{abstract}

\maketitle

\tableofcontents

\section{Introduction}
Group representation theory has been crucial in so many areas of
mathematics that there is essentially no need to speak further of
its successes.  The same is not the case for semigroup
representation theory at the present time.  This is beginning to
change, in a large part due to work of Brown~\cite{Brown1,Brown2}
and Bidigare \textit{et al.}~\cite{BHR} who found applications to
random walks and connections with Solomon's descent algebra.
Solomon, himself, has also found interest in representations of
semigroups, in particular inverse semigroups, for the purposes of
algebraic combinatorics and representation theory of the symmetric
group~\cite{Solomon3,Solomon}. Putcha has applied semigroup
representation theory to finding weights for finite groups of Lie
type~\cite{Putchaweights} and has explored other connections with
modern representation theory~\cite{Putchachar,quiver}. Recently,
Aguiar and Rosas~\cite{Hopf} have used the inverse monoid of uniform
block permutations to study Malvenuto and Reutenauer's Hopf algebra
of permutations and the Hopf algebra of non-commutative symmetric
functions.  Representations of infinite inverse semigroups on
Hilbert spaces have received a lot of attention in the $C^*$-algebra
community, see the book of Paterson~\cite{Paterson} for more
details.  This author has been trying over the past couple of years
to apply semigroup representations to finite semigroup theory and to
automata theory~\cite{myrad,synch}.

One of the great successes of group representation theory is
character theory.  Thanks to Maschke's theorem and the orthogonality
relations and their consequences, much of group representation
theory boils down to the combinatorics of characters and character
sums~\cite{Serre}. Despite intensive work in the fifties and sixties
on representations of finite
semigroups~\cite{CP,Munn,Munnalg,Munnchar,McAlister,McAlistersurvey,RhodesZalc},
culminating in the reduction of calculating irreducible
representations to group representation theory and combinatorics in
matrix algebras over group rings, very few results on characters and
how to calculate multiplicities of irreducible constituents have
been obtained until now.  This is, of course, complicated by the
fact that semigroup algebras are almost never semisimple.  But even
in cases where they are known to be so, like for inverse semigroups,
no nice combinatorial formulas seem to exist in the literature.

There are three notable exceptions.  First there are Munn's
results on characters of the symmetric inverse monoid (also known as
the rook monoid)~\cite{Munnchar}.  These results were extended by
Solomon, who obtained multiplicity formulas for irreducible
constituents using combinatorics associated to partitions, Ferrer's
diagrams, symmetric functions and symmetric groups~\cite{Solomon}.
Then there is Putcha's work on the characters of the full
transformation semigroup~\cite{Putchachar} and his work on monoid
quivers~\cite{quiver}, in which he develops multiplicity formulas
for certain representations. In particular he obtained a formula in
terms of the M\"obius function on the $\J$-order for multiplicities
of irreducible constituents for representations of idempotent
semigroups acting on their left ideals. This same formula, in the
special case of a minimal left ideal, was obtained independently by
Brown~\cite{Brown1,Brown2} and made much more famous because the
formulas were applied to random walks on chambers of hyperplane
arrangements and to other Markov chains to obtain absolutely amazing
results! Also Brown's work, based on the work of Bidigare
\textit{et al.}~\cite{BHR} for hyperplane face semigroups, developed
the theory from scratch, making it  accessible to the general
public.  Putcha's work, however goes much deeper from the
representation theoretic point-of-view: he shows that regular
semigroups have quasi-hereditary algebras and he calculates the
blocks in terms of character formulas~\cite{quiver}.

In our previous paper~\cite{mobius1}, we showed how
Solomon's~\cite{Solomon1} approach to the semigroup algebra of a
semilattice, that is an idempotent inverse semigroup, via the
M\"obius algebra can be extended to inverse semigroups via the
groupoid algebra. This allowed us to obtain an explicit
decomposition of the algebra of an inverse semigroup into a direct
sum of matrix algebras over algebras of maximal subgroups.  Also we
were able to explicitly determine the central primitive idempotents
of the semigroup algebra in terms of character sums and the M\"obius
function of the inverse semigroup.

In this paper we use the above decomposition to give a character
formula for multiplicities of irreducible constituents in
representations of inverse semigroups.  In particular, we recover
and greatly generalize Solomon's results~\cite{Solomon} for the
symmetric inverse monoid concerning decompositions of tensor
and exterior powers of rook matrix representations.  Moreover, we
obtain the results in a more elementary fashion.

 We also give character theoretic proofs of the description of the
decomposition of partial permutation representations into
irreducible constituents (this last result can also be obtained by
the classical semigroup techniques~\cite{CP,RhodesZalc}, with
greater effort).

Just as the irreducible representations of idempotent semigroups
correspond to irreducible representations of semilattices (this has
been known to semigroup theorists
since~\cite{CP,Munnalg,Munn,RhodesZalc} and has recently been
popularized by Brown~\cite{Brown2}), there is a large class of
semigroups whose irreducible representations essentially factor
through inverse semigroups; this includes all finite semigroups with
basic algebras. Our results therefore extend to this domain, and in
particular we recover the case of idempotent
semigroups~\cite{Brown1,Brown2,quiver} and semigroups with basic
algebras~\cite{mobius1}.  In the process we calculate an explicit
basis for the radical of the semigroup algebra of a finite semigroup
with commuting idempotents as well as identifying the semisimple
quotient as a certain retract.

  Our aim is to make this paper accessible to
people interested in algebraic combinatorics, semigroups and
representation theory and so we shall try to keep specialized
semigroup notions to a minimum.  In particular, we shall try to
prove most results from~\cite{CP} that we need, or refer
to~\cite{myrad}, where many results that we need are proved in a
less semigroup theoretic language than~\cite{CP}.
Putcha~\cite{Putchachar} gives a nice survey of semigroup
representation theory, but for the inverse semigroup case our
methods handle things from scratch.

The paper is organized as follows. We begin with a brief
introduction to inverse semigroups.  This is followed by a review of
Rota's theory of incidence algebras and M\"obius inversion.  The
results of our first paper~\cite{mobius1} are then summarized.  The
main argument of~\cite{mobius1} is proved in a simpler (and at the
same time more complete) manner.
 The following section gives the general formula for inverse
semigroup intertwining numbers.  To demonstrate the versatility of
our formula, we compute several examples involving tensor and
exterior products of rook matrix representations.  We then compare
our method for computing multiplicities to Solomon's method
(properly generalized) via character tables. Finally we explain how
to use the inverse semigroup results to handle more general
semigroups. This last section will be more demanding of the reader
in terms of semigroup theoretic background, but most of the
necessary background can be found in~\cite{CP,Arbib,myrad,qtheor}.  In this
last section, we also finish the work begun in~\cite{mobius1} on
analyzing random walks on triangularizable semigroups.  In that
paper, we calculated the eigenvalues, but were unable to determine
multiplicities under the most general assumptions.  In this paper we
can handle the general case.

Herein we adopt the convention that all transformation groups
and semigroups act on the right of sets.  We also consider only
right modules.

\section{Inverse Semigroups}
Inverse semigroups capture partial symmetry in much the same way
that groups capture symmetry; see Lawson's book~\cite{Lawson} for
more on this viewpoint and the abstract theory of inverse
semigroups.

\subsection{Definition and basic properties}
Let $X$ be a set. We shall denote by $\mathfrak S_X$ the symmetric
group on $X$.   If $n$ is a natural number, we shall set
$[n]=\{1,\ldots,n\}$. The symmetric group on $[n]$ will be denoted by
$\mathfrak S_n$, as usual. What is a partial permutation? An example
of a partial permutation of the set $\{1,2,3,4\}$ is \[\sigma =
\begin{pmatrix} 1 & 2 & 3 & 4\\ - & 3 & - &1\end{pmatrix}.\]  The
domain of $\sigma$ is $\{2,4\}$ and the range of $\sigma$ is
$\{1,3\}$.  More formally a \emph{partial permutation} of a set $X$
is a bijection $\sigma:Y\to Z$ with $Y,Z\subseteq X$.  We admit the
possibility that $Y$ and $Z$ are empty.   Partial permutations can
be composed via the usual rule for composition of partial functions
and the monoid of all partial permutations on a set $X$ is called
the \emph{symmetric inverse monoid} on $X$, denoted $\mathfrak I_X$.
The empty partial permutation is the zero element of $\mathfrak
I_X$, and so will be denoted $0$.  We shall write $\mathfrak I_n$
for the symmetric inverse monoid on $[n]$.  Clearly $\mathfrak S_n$
is the group of units of $\mathfrak I_n$.

  For reasons that will be come apparent later, we shall use the term \emph{subgroup}
to mean any subsemigroup of a semigroup that happens to be a group.
For instance, if $Y\subseteq X$, then the collection of all partial
permutations of $X$ with domain and range $Y$ is a subgroup of
$\mathfrak I_X$, which is isomorphic to $\mathfrak S_Y$. This is a
nice feature of $\mathfrak I_n$:  it contains in a natural way all
symmetric groups of degree at most $n$ and its representations
relate to the representations of each of these symmetric groups.

The monoid $\mathfrak I_X$ comes equipped with a natural involution
that takes a bijection $\sigma:Y\to Z$ to its inverse $\sigma\inv:
Z\to Y$. The key properties of the involution are, for
$\sigma,\tau\in \mathfrak I_X$:
\begin{itemize}
\item $\sigma\sigma\inv \sigma = \sigma$;
\item $\sigma\inv \sigma\sigma\inv =\sigma\inv$;
\item $(\sigma\inv)\inv =\sigma$;
\item $(\sigma\tau)\inv = \tau\inv \sigma\inv$;
\item $\sigma\sigma\inv \tau\tau\inv = \tau\tau\inv \sigma \sigma\inv$.
\end{itemize}

Let us define a (concrete) \emph{inverse semigroup} to be an
involution-closed subsemigroup of some $\mathfrak I_X$.   There is
an abstract characterization, due to Preston and
Vagner~\cite{CP,Lawson}, which says that $S$ is an inverse semigroup
if and only if, for all $s\in S$, there is a unique $t\in T$ such
that $sts=s$ and $tst=t$; one calls $t$ the \emph{inverse} of $s$
and denotes it by $s\inv$.  Moreover, $S$ in this case has a
faithful representation (called the Preston-Vagner representation)
\mbox{$\rho_S:S\to \mathfrak I_S$} by partial permutations  via
partial right multiplication~\cite{CP,Lawson}. Thus we can view
finite inverse semigroups as involution closed subsemigroups of
$\mathfrak I_n$ and we shall draw our intuition from there. Later we
shall give a new interpretation of the Preston-Vagner representation
for finite inverse semigroups in terms of their semigroup algebras.

Another way to think of inverse semigroups is via so-called rook
matrices.  An $n\times n$ \emph{rook matrix} is a matrix of zeroes
and ones with the constraint that if we view the ones as rooks on an
$n\times n$ chessboard, then the rooks must all be non-attacking. In
other words a rook matrix is what you obtain from a permutation
matrix by replacing some of the ones by zeroes.  One might equally
well call these partial permutation matrices and it is clear that
the monoid of $n\times n$ rook matrices $R_n$, called by Solomon the
\emph{rook monoid}~\cite{Solomon}, is isomorphic to the symmetric inverse
monoid. If $\sigma\in \mathfrak I_n$, the corresponding element of
$R_n$ has a one in position $i,j$ if $i\sigma =j$ and zero
otherwise. In $R_n$ the involution is given by the transpose.  It
follows that any representation of an inverse semigroup $S$ by rook
matrices is completely reducible over any subfield of $\mathbb{C}$
since if one uses the usual inner product, then the orthogonal
complement of an $S$-invariant subspace will be $S$-invariant (as is
the case for permutation representations of groups).   If $\sigma\in
\mathfrak I_n$, then by the \emph{rank} of $\sigma$, denoted
$\mathrm{rk}(\sigma)$, we mean the cardinality of the range of
$\sigma$; this is precisely the rank of the associated rook matrix.

Let $M_n(K)$ be the monoid of $n\times n$ matrices over a field $K$.
Let $B$ be the Borel subgroup of invertible $n\times n$ upper
triangular matrices over $K$.  Then one easily verifies that the
Bruhat decomposition of $Gl_n(K)$ extends to $M_n(K)$ via $R_n$:
that is $M_n(K) = \biguplus_{r\in R_n} BrB$.  Renner showed
more generally that, for any reductive algebraic monoid, there is a
Bruhat decomposition involving the Borel subgroup of the reductive
group of units and a finite inverse monoid $R$, now called the
Renner monoid~\cite{Putcha,Renner,Rennerbook}.   Solomon defined a Hecke algebra
in this context~\cite{Solomon3} and Putcha~\cite{PutchaHecke}
studied the relationship of the Hecke algebra with the algebra of
the Renner monoid.

\subsection{Idempotents and order}
If $S$ is a semigroup, then we denote by $E(S)$  the set of
idempotents of $S$. Observe that the idempotents of $\mathfrak I_n$
are the partial identities $1_X$, $X\subseteq [n]$. Also the
equality
\begin{equation}\label{islattice}
1_X1_Y = 1_{X\cap Y} = 1_Y1_X
\end{equation}
holds.   Thus $E(\mathfrak I_n)$ is a commutative semigroup
isomorphic to the Boolean lattice $\mathfrak B_n$ of subsets of
$[n]$ with the meet operation. In general, if $S\leq \mathfrak I_n$
is an inverse semigroup, then $E(S)$ is a meet subsemilattice of
$\mathfrak B_n$ (it will be a sublattice if $S$ is a submonoid). The
order can be defined intrinsically by observing that
\begin{equation}\label{idemorder}
e\leq f \iff ef=e=fe,\ \text{for}\ e,f\in E(S).
\end{equation}

The ordering on idempotents extends naturally to the whole
semigroup: let us again take $\mathfrak I_n$ as our model.  We can
define $\sigma\leq \tau$, for $\sigma,\tau\in \mathfrak I_n$, if
$\sigma$ is a restriction of $\tau$.  This order is clearly
compatible with multiplication and $\sigma\leq \tau$ implies
$\sigma\inv \leq \tau\inv$.  Thus if $S\leq \mathfrak I_n$ is an
inverse semigroup, it too has an ordering by restriction. Again the
order can be defined intrinsically:  it is easy to see that if
$s,t\in S$, then
\begin{equation}\label{order}
s\leq t\iff s = et,\ \text{some}\ e\in E(S)\iff s =tf,\ \text{some}\
f\in E(S).
\end{equation}
If one likes, one can take~\eqref{order} as the definition of the
\emph{natural partial order} on $S$~\cite{CP,Lawson}.

Let us use the notation $\mathrm{dom}(\sigma)$ for the domain of
$\sigma\in \mathfrak I_n$ and $\mathrm{ran}(\sigma)$ for the range.
If $\sigma \in \mathfrak I_n$, then (recalling that $\mathfrak I_n$
acts on the right of $[n]$)
\begin{align}\label{domain}
\sigma\sigma\inv &=1_{\mathrm{dom}(\sigma)}\\
\label{range} \sigma\inv \sigma & = 1_{\mathrm{ran}(\sigma)}.
\end{align} Thus if $s$ is an element of an inverse semigroup, then it is natural to think
of $ss\inv$ as the ``domain" of $s$ and $s\inv s$ as the ``range" of
$s$ and so we shall write
\begin{equation}\label{domainsandrange}
\begin{split}
ss\inv &= \dom s\\
s\inv s &=\ran s
\end{split}
\end{equation}
This means that we are going to abuse the distinction between a
partial identity and the corresponding subset.  So if $S\leq
\mathfrak I_n$ and $s\in S$, we shall write \mbox{$x\in \dom s$} to mean $x$
belongs to the domain of $s$.  With this viewpoint it is natural to
think of $s$ as an isomorphism from $\dom s$ to $\ran s$.  So let us
define \mbox{$e,f\in E(S)$} to be \emph{isomorphic} if there exists
$s\in S$ with $\dom s = e$ and $\ran s =f$: that is $e=ss\inv$,
$f=s\inv s$. Following long standing semigroup tradition, going back
to Green~\cite{Green,CP,Lawson}, we shall write $e\D f$.  One can
extend this relation to all of $S$ by defining $s\D t$ if $\dom s$
is isomorphic to $\dom t$ (or equivalently $\ran s\D \ran t$).  One
can easily verify that $\D$ is an equivalence relation on $S$. With
a bit more work, one can verify that if $S$ is a finite inverse
semigroup, then $s\D t$ if and only if $s$ and $t$ generate the same
two-sided ideal~\cite{CP,Lawson}. The equivalence classes with
respect to the $\D$-relation are called \emph{$\D$-classes},
although connected components in this context would be a better
word.

Let $S\leq\mathfrak I_n$ and let $e\in E(S)$.  Then $e$ is the
identity of a subset $X\subseteq [n]$.  It is then clear that the
set
\[G_e = \{s\in S\mid \dom s = e = \ran s\}\] is a permutation group
of degree $\rk e$.  Actually $G_e$ makes perfectly good sense via
\eqref{domainsandrange}, and is a group, without any reference to an
embedding of $S$ into $\mathfrak I_n$. It is called the
\emph{maximal subgroup} of $S$ at $e$, as it contains any other subgroup of $S$ with identity $e$.  It is straightforward to see
that if $e,f\in E(S)$ are isomorphic idempotents, then $G_e\cong
G_f$.  In fact if $s\in S$ with \mbox{$\dom s= e,\ \ran s = f$},
then conjugation by $s$ implements the isomorphism. It was first
observed by Munn and Ponizovski{\u\i}~\cite{Munnalg,Munn,CP,Poni} that the representation
theory of $S$ is in fact controlled by the representations of its
maximal subgroups. We shall see this more explicitly below.

Recall that an \emph{order ideal} in a partially ordered set $P$ is
a subset $I$ such that $x\leq y\in I$ implies $x\in I$.  If $p\in
P$, then $p^{\downarrow}$ denotes the \emph{principal order ideal}
generated by $p$.  So \[p^{\downarrow} = \{x\in P\mid x\leq p\}.\]
As usual, for $p_1,p_2\in P$, the closed interval from $p_1$ to
$p_2$ will be denoted $[p_1,p_2]$.

It is easy to see that if $S$ is an inverse semigroup then the
idempotent set $E(S)$ is an order ideal of $S$: any restriction of a
partial identity is a partial identity. So if $e,f\in E(S)$, then
the interval $[e,f]$ in $S$ and in $E(S)$ coincide.  Also one can
verify that if $s,t\in S$, then the following intervals are
isomorphic posets:
\begin{equation}\label{sameintervals}
[s,t]\cong [\dom s,\dom t]\cong [\ran s,\ran t].
\end{equation}

\section{Incidence Algebras and M\"obius Functions}
Let $(P,\leq)$ be a finite partially ordered set and $A$ a
commutative ring with unit.   The \emph{incidence algebra} of $P$
over $A$, which we denote $A\pow {P}$, is the algebra of all
functions $f:P\times P\to A$ such that \[f(x,y)\neq 0 \implies x\leq
y\] equipped with the convolution product
\[(f\ast g)(x,y) = \sum_{x\leq z\leq y}f(x,z)g(z,y).\]  In other
words one can think of $A\pow {P}$ as the algebra of all $P\times P$
upper triangular matrices over $A$, where upper triangular is
defined relative to the partial order on $P$.

With this product and pointwise addition $A\pow {P}$ is an
$A$-algebra with unit the Kronecker delta function
$\delta$~\cite{Hall,Stanley}. An element $f\in A\pow {P}$ is
invertible if and only if $f(x,x)$ is a unit of $A$ for all $x\in
P$~\cite{Hall,Stanley}.  One can define the inverse inductively by
\begin{equation}\label{incidenceinverse}
\begin{split}
f\inv (x,x) &= f(x,x)\inv\\
f\inv (x,y) &= -f(x,x)\inv \sum_{x<z\leq y} f(x,z)f\inv (z,y),\
\text{for}\ x\leq y.
\end{split}
\end{equation}

The \emph{zeta function} $\zeta$ of $P$ is the element of $A\pow
{P}$ that takes on the value of $1$ whenever $x\leq y$ and $0$
otherwise. The zeta function is invertible over any ring $A$ and its inverse is
called the \emph{M\"obius function}.  The M\"obius function only depends on the
characteristic; the characteristic zero version is called the
\emph{M\"obius function} of $P$ and is denoted $\mu$, or $\mu_P$ if
we wish to emphasize the partially ordered set.
   From~\eqref{incidenceinverse}, it follows that $\mu(x,y)$
depends only on the isomorphism class of the interval $[x,y]$.  In
particular, for an inverse semigroup $S$, the M\"obius function for
$S$ is determined by the M\"obius function for $E(S)$
via~\eqref{sameintervals}. The following is Rota's M\"obius
Inversion Theorem~\cite{Hall,Stanley}.

\begin{Thm}[M\"obius Inversion Theorem]\label{inversion}
Let $(P,\leq)$ be a finite partially ordered set  and $G$ be an
Abelian group. Suppose that $f:P\to G$ is a function and define
$g:P\to G$ by
\[g(x) = \sum_{y\leq x} f(y).\] Then \[f(x) = \sum_{y\leq x} g(y)\mu
(y,x).\]
\end{Thm}

Returning to our motivating example $\mathfrak I_n$, the M\"obius
function for $\mathfrak B_n$ is well known~\cite{Hall,Stanley}: if
$Y\subseteq Z$, then
\[\mu_{\mathfrak B_n}(Y,Z) = (-1)^{|Z|-|Y|}.\]  Hence, for $\mathfrak I_n$,
the M\"obius function is determined by
\begin{equation}\label{symmetric}
\mu_{\mathfrak I_n}(s,t) = (-1)^{\rk t-\rk s}
\end{equation}
for $s\leq t$.

The semilattice of idempotents of any Renner monoid $R$ is the face
lattice of a rational polytope~\cite{Putcha,Rennerbook}.  Such a
partially ordered set is Eulerian~\cite{Stanley} and so has a
well-defined rank function. This rank function extends to $R$ and in
this situation the M\"obius function is given by~\eqref{symmetric}.

\section{Inverse Semigroup Algebras}
Let $A$ be a commutative ring with unit.  Solomon~\cite{Solomon1}
assigned to each partially ordered set $P$ an $A$-algebra called the
M\"obius algebra of $P$. When $P$ is a semilattice, he used M\"obius
inversion to show that the M\"obius algebra is isomorphic to the
semigroup algebra of $P$. In~\cite{mobius1} we extended this to
arbitrary finite inverse semigroups. We review the construction here
from a different viewpoint; further details can be found
in~\cite{mobius1}.

\subsection{The groupoid algebra}
Let $S$ be a finite inverse semigroup.  We define an $A$-algebra
called the groupoid algebra of $S$ over $A$ for reasons that will become
clear.  Let us denote by $G(S)$ the set $\{\lfloor s \rfloor\mid
s\in S\}$, a formal disjoint copy of $S$. To motivate the definition
of the groupoid algebra, let us point out that there is another way
to model partial bijections: allowing composition if and only if the
domains and ranges line up.  The groupoid algebra encodes this.  So
let $AG(S)$ be a free $A$-module with basis $G(S)$ and define a
multiplication on $AG(S)$ by setting, for $s,t\in S$,
\begin{equation}\label{groupoidmult}
\lfloor s \rfloor\lfloor t \rfloor = \begin{cases} \lfloor st
\rfloor & \text{if}\ \ran s=\dom t
\\ 0 & \text{else.}\end{cases}
\end{equation}
It is easy to check that~\eqref{groupoidmult} extends to an
associative multiplication on $AG(S)$.   In fact, the set $G(S)$
with the above multiplication (where $0$ is interpreted as
undefined) is a groupoid in the sense of a small
category~\cite{MacLane} in which all arrows are invertible;
see~\cite[Section 3.1]{Lawson} for details.  We shall use this
fact and its consequences without further comment. In particular, we
use that if $\lfloor s \rfloor\lfloor t \rfloor$ is defined, then
$\dom {st} = \dom {s}$ and $\ran {st} = \ran {t}$. The algebra
$AG(S)$ is termed the \emph{groupoid algebra} of
$S$~\cite{mobius1}.  Notice that $\sum_{e\in E(S)} \lfloor e\rfloor$ is an identity for $AG(S)$ and so $AG(S)$ is unital.

Define temporarily, for $s\in S$, an element $v_s\in AG(S)$ by \[v_s
= \sum_{t\leq s} \lfloor t \rfloor.\]  So $v_s$ encodes $s$ via its
restrictions.  Then by M\"obius inversion:
\[\lfloor s \rfloor = \sum_{t\leq s} v_t\mu(t,s).\]  Hence
$\{v_s\mid s\in S\}$ is a basis for $AG(S)$.  The following result
from~\cite{mobius1} can be viewed as a coordinate-free form of the
Sch\"utzenberger representation~\cite{Schutzrep} or Preston-Vagner
representation~\cite{CP,Lawson}. We include the proof for
completeness. The proof here includes a small detail that was
unfortunately omitted from~\cite{mobius1}. Nonetheless, the proof is
more compact, as we avoid considering explicitly the
direct sum of Sch\"utzenberger representations, as was done
in~\cite{mobius1}.

\begin{Lemma}\label{groupoid=semigroup}
Let $s,t\in S$. Then $v_sv_t = v_{st}$.
\end{Lemma}
\begin{proof}
First we compute
\begin{align*}
v_sv_t &= \left(\sum_{s'\leq s} \lfloor s'
\rfloor\right)\left(\sum_{t'\leq t} \lfloor t' \rfloor\right)
\\ &= \sum_{s'\leq s,\ t'\leq t,\ \ran {s'}=\dom
{t'}} \lfloor s't' \rfloor.
\end{align*}
Since the natural partial order is compatible with multiplication if
$s'\leq s$ and $t'\leq t$, it follows that $s't'\leq st$. Thus to
obtain the desired result it suffices to show that each $u\leq st$
can be  uniquely factored as a product $s't'$ with $s'\leq s$, $t'\leq
t$ and $\ran{s'}=\dom{t'}$.  Note that $u\leq st$ implies that
$uu\inv st = u = stu\inv u$.

To obtain such a factorization $u=s't'$ we must have $\dom{s'}=\dom
u$ and $\ran {t'}=\ran u$. That is, we must have $s'=uu\inv s$ and
$t'=tu\inv u$.  Clearly $s't' = uu\inv stu\inv u = uu\inv u =u$.
Let us check that $\ran {s'}=\dom {t'}$. First observe that $(s')\inv
= tu\inv$. Indeed,
\begin{align*}
s'(tu\inv) s' &= (uu\inv stu\inv u)u\inv s = uu\inv s=s'
\\
(tu\inv) s'(tu\inv) &= tu\inv (uu\inv st)u\inv = tu\inv uu\inv =
tu\inv.
\end{align*}
  Similarly,
one can verify that $(t')\inv = u\inv s$.
Thus we see \[(s')\inv s' = tu\inv (uu\inv s) = (tu\inv u)u\inv s=
t'(t')\inv,\] as desired.
\end{proof}

As a consequence, it follows that $AG(S)\cong AS$ and in particular $AS$ is unital.  Henceforth, we shall
identify $AS$ with $AG(S)$ by identifying $s$ with $v_s$ and so from now on the notation $AG(S)$ will be dropped. Our viewpoint is
that the semigroup algebra $AS$ has two natural bases: the
usual basis and the basis $\{\lfloor s \rfloor\mid s\in S\}$. Let us
make this precise~\cite{mobius1}.

\begin{Thm}\label{mainthm:mobius}
Let $A$ be a unital commutative ring and $S$ a finite inverse semigroup.
Define, for $s\in S$,
\begin{equation}\label{definenewbasis}
\lfloor s \rfloor = \sum_{t\leq s} t\mu(t,s).
\end{equation}
Then $\{\lfloor s \rfloor\mid s\in S\}$ is a basis for $AS$ and the
multiplication with respect to this basis is given
by~\eqref{groupoidmult}.
\end{Thm}

It is worth remarking that Theorem~\ref{mainthm:mobius} remains
valid for infinite inverse semigroups $S$ so long as one assumes
descending chain condition on the set of idempotents of $S$.
  We now give a new proof of the semisimplicity of $KS$
for $S$ a finite semigroup and $K$ a subfield of $\mathbb{C}$ by
comparing the linear representations associated to the two bases of
$KS$.  If we use the usual basis $\{s\mid s\in S\}$, then the
associated linear representation is just the regular representation
of $S$. What we now wish to show is that if we use the basis
$\{\lfloor s\rfloor\mid s\in S\}$, then we obtain the rook matrix
representation associated to the Preston-Vagner representation
$\rho_S:S\to \mathfrak I_S$. In particular, these two
representations will be equivalent and so, since the latter (as
observed earlier) is completely reducible, we shall obtain that the
regular representation of $KS$ is completely reducible, whence $KS$
is a semisimple algebra.

Let us recall the definition of $\rho_S:S\to \mathfrak I_S$; details
can be found in~\cite{CP,Lawson,Paterson}. For $s,t\in S$,
\[s\rho_S(t) =
\begin{cases} st & s\inv s\leq tt\inv \\ \text{undefined} &
\text{else.}\end{cases}\]  Our claim in the above paragraph is then
an immediate consequence of the following proposition.

\begin{Prop}
Let $S$ be a finite inverse semigroup and $A$ a commutative ring
with unit.  Then for $s,t\in S$, \[\lfloor s\rfloor t =\begin{cases}
\lfloor st\rfloor & s\inv s\leq tt\inv \\ 0 &
\text{else.}\end{cases}\]
\end{Prop}
\begin{proof}
By M\"obius inversion, $t= \sum_{u\leq t}\lfloor u\rfloor$.  So
\begin{align*}
\lfloor s\rfloor t & = \lfloor s\rfloor \sum_{u\leq t}\lfloor
u\rfloor\\ &=\sum_{u\leq t}\lfloor s\rfloor \lfloor u\rfloor.
\end{align*}
 Now $\lfloor s\rfloor \lfloor u\rfloor = \lfloor su\rfloor$ if
$s\inv s = uu\inv$ and is zero else.  There can be at most one
element $u\leq t$ with $s\inv s = uu\inv$~\cite[Theorem 3.1.2]{Lawson}.
 Now if $u\leq t$ and $s\inv s = uu\inv$, then
$s\inv s\leq tt\inv$. Conversely, if $s\inv s\leq tt\inv$, define
$u=s\inv st$. Then $u\leq t$ and
\[uu\inv = s\inv stt\inv s\inv s = s\inv s.\]
Moreover, $su = s(s\inv st) = st$.  Thus
\[\lfloor s\rfloor t=\sum_{u\leq t}\lfloor s\rfloor \lfloor u\rfloor
= \begin{cases} \lfloor st\rfloor & s\inv s\leq tt\inv \\
0&\text{else,}\end{cases}\] completing the proof of the proposition.
\end{proof}

\begin{Cor}
Let $S$ be a finite inverse semigroup and $K$ any field.  Then the
regular representation of $S$ and the Preston-Vagner representation
of $S$ are equivalent as linear representations.  In particular, if
$K$ is a subfield of $\mathbb{C}$, then $KS$ is semisimple.
\end{Cor}

\subsection{Decomposition into matrix algebras over group rings}
We now turn to decomposing $AS$ into matrix algebras over group
algebras. Let $S$ be a finite inverse semigroup with $\D$-classes
$D_1,\ldots, D_r$. Recall these are the equivalence classes
corresponding to isomorphic idempotents. Let $AD_i$ be the $A$-span
of $\{\lfloor s \rfloor\mid s\in D_i\}$. The following result is
immediate from~\eqref{groupoidmult} and
Theorem~\ref{mainthm:mobius}.

\begin{Thm}\label{decomp1}
Let $A$ be a unital commutative ring and let $S$ be a finite inverse
semigroup with $\D$-classes $D_1,\ldots,D_r$.  Then $AS =
\bigoplus_{i=1}^r AD_i$.
\end{Thm}

For each $i$, fix an idempotent $e_i$ of $D_i$ and let $G_{e_i}$ be
the corresponding maximal subgroup.  Since the idempotents of $D_i$
are isomorphic, this group does not depend on the choice of $e_i$ up
to isomorphism.  Let $n_i = |E(D_i)|$, that is $n_i$ denotes the
number of idempotents in $D_i$.  The structure of the algebra $AD_i$ can best be understood via the following simple argument.  As a right $AD_i$-module, we have $AD_i\cong n_ie_iAD_i$ and hence \[AD_i\cong \mathrm{End}_{AD_i}(n_ie_iAD_i) \cong M_{n_i}(e_iAD_ie_i) \cong M_{n_i}(AG_{e_i}).\] To make the isomorphism more explicit, we recall the argument
from~\cite[Theorem 3.2]{mobius1} (essentially due to
Munn and Ponizovski{\u\i}~\cite{Munnalg,Munn,CP,Poni}) that \mbox{$AD_i\cong
M_{n_i}(AG_{e_i})$}. We view $n_i\times n_i$ matrices as having rows and columns indexed by
the idempotents of $D_i$. Fix, for each $e\in D_i$, an element
$p_e\in S$ with $\dom {p_e} = e_i$ and $\ran {p_e} =e$. We take
$p_{e_i}=e_i$. Define a map $\p:AD_i\to M_{n_i}(AG_{e_i})$ on a
basis element $\lfloor s \rfloor\in AD_i$ with $\dom s =e$ and $\ran
s=f$ by
\begin{equation}\label{fundgroup}
\p(\lfloor s \rfloor) = p_{e}sp_{f}\inv E_{e,f}
\end{equation}
 where $E_{e,f}$ is
the standard matrix unit with $1$ in position $e,f$ and zero in all
other positions. Observe that $p_{e}sp_{f}\inv\in G_{e_i}$ by
construction. It is straightforward~\cite{mobius1} to show that $\p$
is an isomorphism and to verify that the inverse is induced by
$gE_{ef}\mapsto \lfloor p_e\inv gp_f\rfloor$, for $g\in G_{e_i}$ and $e,f$
idempotents in $D_i$. The reader should compare~\eqref{fundgroup} to
the calculation of the fundamental group of a graph.

As a consequence of the above isomorphism, we obtain the following
result, which is implicit in the work of Munn and Ponizovski{\u\i}~\cite{Munn,Munnalg,CP,Poni}
and can be found explicitly in~\cite{Okninski} and~\cite[Theorem 4.3]{mobius1}.

\begin{Thm}\label{matrixdecomp}
Let $S$ be a finite inverse semigroup with representatives
$e_1,\ldots, e_r$ of the isomorphism classes of idempotents.  Let
$n_i$ be the number of idempotents isomorphic to $e_i$.  Let $A$ be
a commutative ring with unit.  Then $AS\cong \bigoplus_{i=1}^r
M_{n_i}(AG_{e_i})$.
\end{Thm}

This decomposition implies the well known fact that the size of $S$
is $\sum_{i=1}^r n_i^2|G_{e_i}|$.   One may also deduce the
following theorem~\cite{Munn,Poni,CP}.

\begin{Cor}[Munn, Ponizovski{\u\i}]
Let $K$ be a field and $S$ a finite inverse semigroup.  If the
characteristic of $K$ is $0$, then $KS$ is semisimple.  If the
characteristic of $K$ is a prime $p$, then $KS$ is semisimple if and
only if $p$ does not divide the order of any maximal subgroup of
$S$.
\end{Cor}

We mention that Solomon~\cite{Solomon} gives the exact decomposition
obtained above for the special case of $\mathfrak I_n$ (this
decomposition was obtained independently by V.~Dlab in unpublished
work). In~\cite{mobius1} we went on to describe explicitly central
idempotents and central primitive idempotents. However, the author,
somewhat embarrassingly, missed that the above decomposition allows
one to obtain a formula for character multiplicities. (Solomon also
seems to have missed this~\cite{Solomon} since he uses a different
approach to obtain multiplicity formulas for $\mathfrak I_n$; we
shall compare the two approaches in Section~\ref{compare}.) The goal
of the next few sections is to rectify this.

\section{Character Formulas for Multiplicities}
In this section we assume that $K$ is a field of characteristic
zero.   The most interesting case is when $K=\mathbb{C}$, the
complex field.  If $\p$ is a representation of a group, semigroup or algebra, then the \emph{character} of $\p$ is the composition $\mathrm{tr}\circ \p$ where $\mathrm{tr}$ is the trace.
Fix a finite inverse semigroup $S$ and, for each
$\D$-class $D_1,\ldots, D_r$ of $S$, fix an idempotent $e_i$ and set
$G_i=G_{e_i}$.  Again, let $n_i$ be the number of idempotents in
$D_i$.

It is clear from Theorem~\ref{matrixdecomp} that the algebra $KS$ is Morita equivalent to \mbox{$KG_1\times \cdots\times KG_r$} and so the irreducible
representations of $S$ over $K$ correspond to elements of the set
$\biguplus_{i=1}^r \mathrm{Irr}(G_i)$, where $\mathrm{Irr}(G_i)$ is
the set of irreducible representations of $G_i$ (up to equivalence).
Namely, if $\p$ is an irreducible representation of $G_i$, then we
can tensor it up to $KD_i\cong M_{n_i}(KG_i)$ and then extend to
$KS$ by making it zero on the other summands.  Let us call this
representation $\p^*$.  If $\chi$ is the character of $G$ associated
to $\p$, denote by $\chi^*$ the character of $S$ associated to
$\p^*$.  Note that $\mathrm{deg}(\p^*) = n_i\mathrm{deg}(\p)$. Thus
over an algebraically closed field we see that, as in the group
case, the size of $S$ is the sum of the squares of the degrees of
the (inequivalent) irreducible representations of $S$.

 Let $G$ be a finite group and
$\psi,\alpha:G\to K$ be functions. We consider the usual bilinear form on the space of
$K$-valued functions on $G$ given by
\[(\psi,\alpha)_G = \frac{1}{|G|}\sum_{g\in
G}\psi(g\inv)\alpha(g).\] If $\chi$ is an irreducible character of
$G$ and $\alpha$ is any character of $G$, then standard group
representation theory says that $(\chi,\alpha)_G = md$ where $m$ is
the multiplicity of $\chi$ as a constituent of $\alpha$ and $d$ is
the degree over $K$ of the division algebra of $KG$-endomorphisms of
the simple $KG$-module corresponding to $\chi$.

Let us consider now an irreducible representation of $M_n(KG)$.  We
identify $G$ with $G\cdot E_{1,1}$. If we have a character $\theta$
of $M_n(KG)$ and an irreducible character $\chi$ of $G$, then
clearly
\begin{equation}\label{matrixform}
(\chi,\theta|_G)_G=md
\end{equation}
 where $m$ is the multiplicity of
$\chi^*$ in $\theta$ and $d$ is the dimension of the associated
division algebra to $\chi$.

We now wish to generalize the above formula to our finite inverse
semigroup $S$.   Given a character $\theta$ of $S$ and an idempotent
$f\in S$, define $\theta_f$ by $\theta_f(s) = \theta(fs)$. If $e\in
E(S)$, then $\theta_f$ restricts to a $K$-valued function on $G_e$
(for which we use the same notation); it need not in general be a
class function. Now if $\chi$ is a $K$-valued function on $G_e$,
define
\begin{equation}\label{semigroupintertwine}
(\chi,\theta)_S = \sum_{f\leq e}(\chi,\theta_f)_{G_e}\mu(f,e).
\end{equation}
Of course, if $S$ is a group this reduces to the usual formula. Here
is the main new result of this paper.

\begin{Thm}\label{multformula}
Let $S$ be a finite inverse semigroup and $e\in E(S)$.  Let $\theta$
be a character of $S$ and let $\chi$ be an irreducible character of
$G_e$.   Denote by $d$ the dimension of the division algebra
associated to $\chi$ and by $m$ the multiplicity of the induced
irreducible character $\chi^*$ of $S$ in $\theta$. Then
\[(\chi,\theta)_S = \sum_{f\leq e}(\chi,\theta_f)_{G_e}\mu(f,e)=md.\]
\end{Thm}
\begin{proof}
Our discussion above about matrix algebras over group algebras and
the explicit isomorphism we have constructed between $KS$ and the
direct sum of matrix algebras over maximal subgroups tells us how to
calculate multiplicities using the basis $\{\lfloor s \rfloor\mid
s\in S\}$ for $KS$. Namely~\eqref{matrixform} is transformed via our
isomorphism to:
\begin{align}
md &= \frac{1}{|G_e|}\sum_{g\in G_e}\chi(g\inv)\theta(\lfloor g \rfloor)\\
 & = \frac{1}{|G_e|}\sum_{g\in G_e}\chi(g\inv)\theta\left(\sum_{t\leq
 g}t\mu(t,g)\right)\\ \label{leftoff}
 & = \frac{1}{|G_e|}\sum_{g\in G_e}\chi(g\inv)\sum_{t\leq
 g}\theta(t)\mu(t,g)
\end{align}
But recall that the order ideal $g^{\downarrow}$ is isomorphic to
the order ideal $e^{\downarrow}$ via the map sending $f\in
e^{\downarrow}$ to $fg\in g^{\downarrow}$.  In particular, for
$f\leq e$, $\mu(fg,g)=\mu(f,e)$.  Thus the right hand side
of~\eqref{leftoff} is equal to
\begin{align*}
\frac{1}{|G_e|}\sum_{g\in G_e}\chi(g\inv)\sum_{f\leq
 e}\theta(fg)\mu(f,e) & = \sum_{f\leq
 e}\left(\frac{1}{|G_e|}\sum_{g\in
 G_e}\chi(g\inv)\theta_f(g)\right)\mu(f,e)\\ & = \sum_{f\leq
 e}(\chi,\theta_f)_{G_e}\mu(f,e) \\ &= (\chi,\theta)_S.
\end{align*}
Thus we obtain $(\chi,\theta)_S = md$, as desired.
\end{proof}

Of course, if $K$ is algebraically closed, or more generally if
$\chi$ is an absolutely irreducible character, then
$(\chi,\theta)_S$ is the multiplicity of $\chi^*$ in $\theta$.  It is well known that for any finite dimensional semisimple algebra over a field of characteristic $0$, a representation is determined by its character~\cite{CR}.  Nonetheless, Theorem~\ref{multformula} yields a direct proof of this for the case of an inverse semigroup.

\begin{Cor}
Let $K$ be a field of characteristic $0$ and let $S$ be a finite
inverse semigroup.  Then two representations of $S$ are equivalent
if and only if they have the same character.
\end{Cor}
\begin{proof}
Since $KS$ is semisimple, two representations are equivalent if and
only if they have the same multiplicity for each irreducible
constituent.  But Theorem~\ref{multformula} shows that the
multiplicities of the irreducible constituents depend only on the
character.
\end{proof}

\section{Applications to Decomposing Representations}
In this section we use the formula from Theorem~\ref{multformula} to
calculate multiplicities in different settings.  Throughout this
section we assume that $K$ is an algebraically closed field of
characteristic zero.

Let us begin with $S=\mathfrak I_n$.  We choose the idempotents $1_{[r]}$ with
$r\leq n$ as a transversal to the set of isomorphism classes of idempotents of $S$.  We identify the maximal subgroup at $1_{[r]}$ with
$\mathfrak S_r$ (where $\mathfrak S_0$ is the trivial group and the
empty set is viewed as having a unique partition).  If $\mu$ is a
partition of $r$, let $\chi^{\mu}$ denote the irreducible character
of $\mathfrak S_r$ associated to $\mu$~\cite{sagan}. Then, for
$\theta$ a character of $\mathfrak I_n$, we obtain the following
multiplicity formula:
\begin{equation}\label{multsymme}
(\chi^{\mu},\theta)_{\mathfrak I_n} = \sum_{X\subseteq
[r]}(-1)^{r-|X|}(\chi^{\mu},\theta_{1_X})_{\mathfrak S_r}
\end{equation}
This should be contrasted with~\cite[Lemma~3.17]{Solomon}, where
Solomon has a different formula that is more difficult to apply.

\subsection{Tensor powers}
Solomon decomposes the tensor powers of the rook matrix
representation of $\mathfrak I_n$ using his formula.  We do the same
using ours, but in a more general setting.  In particular, we can handle wreath
products of the form $G\wr S$ with $G$ a finite group and $S\leq
\mathfrak I_n$ a finite inverse semigroup containing all the
idempotents of $\mathfrak I_n$. This includes, in addition to the
symmetric inverse monoid, the signed symmetric inverse monoid $\ZZ/2\ZZ\wr \mathfrak I_n$.

First we need a standard lemma about Boolean algebras whose proof we
leave to the reader.  We use $x^c$ for the complement of $x$ in a
Boolean algebra.  We continue to use multiplicative notation for the
meet in a semilattice.

\begin{Lemma}\label{product}
Let $B$ be a Boolean algebra and $x\in B$.  Then $B\cong
x^{\downarrow}\times (x^c)^{\downarrow}$ via the maps
\begin{align*}
y&\longmapsto (yx,yx^c)\\
y\vee z&\longmapsfrom (y,z)
\end{align*}
\end{Lemma}

To apply Lemma~\ref{product} we use the well known
fact~\cite{Stanley} that if $P_1,P_2$ are finite partially ordered
sets, then
\begin{equation}\label{mobiusproduct}
\mu_{P_1\times P_2}((p_1,p_2),(p_1',p_2')) =
\mu_{P_1}(p_1,p_1')\mu_{P_2}(p_2,p_2').
\end{equation}
Combining Lemma~\ref{product} and~\eqref{mobiusproduct} we obtain:

\begin{Cor}\label{booleanmob}
Let $B$ be a Boolean algebra and fix $x\in B$.  Then for
\mbox{$a,b\in B$}, $\mu(a,b) = \mu(ax,bx)\mu(ax^c,bx^c)$.
\end{Cor}

In the case of $\mathfrak B_n$, Corollary~\ref{booleanmob} just
asserts that if $X\subseteq [n]$ is fixed,  then for $Y\subseteq Z$,
$(-1)^{|Z|-|Y|} = (-1)^{|Z\cap X|-|Y\cap X|}\cdot (-1)^{|Z\cap
X^c|-|Y\cap X^c|}$.

We shall also need the following reformulation of the fact that
$\mu$ is the inverse of $\zeta$:
\begin{equation}\label{mobiusrestate}
\sum_{x\leq y\leq z}\mu(y,z) = \begin{cases} 1 & x=z\\ 0 &
x<z\end{cases}
\end{equation}

If $s\in S\leq \mathfrak I_n$, then $\mathrm{Fix}(s)$ denotes the
set of fixed points of $s$ on $[n]$.  Of course,
$\mathrm{Fix}(s)\subseteq \dom s$.  Let us establish some terminology.  If $S\leq \mathfrak I_n$, then by the rook matrix representation of $S$, we mean the linear representation associated to the natural partial permutation action of $S$ on $[n]$.  The character of this representation counts the number of fixed points and is hence referred to as the \emph{fixed-point character} of $S$.

\begin{Prop}\label{booleanalgebra}
Let $S\leq \mathfrak I_n$ and $\theta^p$ be the character of the
$p^{th}$-tensor power of the rook matrix representation of $S$. Let
$e\in E(S)$ and $\chi\in \mathrm{Irr}(G_e)$.   Suppose that
$\{X\subseteq \dom e\mid 1_X\in E(S)\}$ is closed under relative complement (i.e.\ $X\mapsto \dom
e\setminus X$). Suppose further that, for all $g\in G_e$,
$1_{\mathrm{Fix}(g)}\in E(S)$. Then
\begin{equation}\label{booleancase}
(\chi,\theta^p)_S = \frac{1}{|G_e|}\mathrm{deg}(\chi)\sum_{f\leq e}
\rk f^p\mu(f,e)
\end{equation}
\end{Prop}
\begin{proof}
We begin by computing:
\begin{equation}\label{leftoffagain}
\begin{split}
(\chi,\theta^p)_S &= \sum_{f\leq
e}(\chi,\theta^p_f)_{G_e}\mu(f,e)\\
&= \sum_{f\leq e}\frac{1}{|G_e|}\sum_{g\in
                  G}\chi(g\inv)\theta^p(fg)\mu(f,e)\\
                  &= \frac{1}{|G_e|}\sum_{g\in
                  G}\chi(g\inv)\sum_{f\leq e}\theta^p(fg)\mu(f,e)
\end{split}
\end{equation}
So let us analyze the term $\sum_{f\leq e}\theta^p(fg)\mu(f,e)$.
Setting $h=1_{\mathrm{Fix}(g)}\in E(S)$ and $h^c=1_{\dom e\setminus
\mathrm{Fix}(g)}\in E(S)$, we can rewrite $\theta^p(fg)$ as follows:
\begin{equation*}\label{eq:tensor}
\theta^p(fg) = |\mathrm{Fix}(fg)|^p = |\mathrm{Fix}(g|_{\dom
f})|^p=|\mathrm{Fix}(g)\cap \dom{f}|^p = \rk {hf}^p.
\end{equation*}
 The above equation, together with
Corollary~\ref{booleanmob} and the fact that $e^{\downarrow}$ is a
Boolean algebra (via our hypotheses), allows us to rewrite our sum:
\begin{align*}
\sum_{f\leq e}\theta^p(fg)\mu(f,e) &= \sum_{x\leq h,\ y\leq
h^c}\rk x^p\mu(x,h)\mu(y,h^c)\\
            &=\sum_{x\leq h}\rk x^p\mu(x,h)\sum_{0\leq y\leq
            h^c}\mu(y,h^c)
\end{align*}
But by~\eqref{mobiusrestate},
\[\sum_{0\leq y\leq
            h^c}\mu(y,h^c)=0\] unless $h^c =0$, or equivalently,
            unless $e=h$. In this latter case,
            we then have that
            $\dom{e}=\mathrm{Fix}(g)$ and hence $g=e$.  Thus the
            right hand side of~\eqref{leftoffagain} becomes:
\begin{equation*} \frac{1}{|G_e|}\chi(e)\sum_{f\leq e}\rk
f^p\mu(f,e) = \frac{1}{|G_e|}\mathrm{deg}(\chi)\sum_{f\leq e}
\rk f^p\mu(f,e).
\end{equation*}
This proves~\eqref{booleancase}
\end{proof}

We now wish to obtain Solomon's result~\cite[Example~3.18]{Solomon} in a
more general form, in particular for wreath products.  Let $G$ be a
finite group and $S\subseteq \mathfrak I_n$.  Their \emph{wreath
product} $G\wr S$ is an inverse semigroup of partial permutations of
$G\times [n]$. It consists of all partial permutations of $G\times
[n]$ that can be expressed in the form $(f,\sigma)$ where $\sigma\in
S$, $f:[n]\to G$ and
\[(g,i)(f,\sigma) = \begin{cases} (gf(i),i\sigma) & \text{if}\
i\sigma\ \text{is defined}\\ \text{undefined} &
\text{else.}\end{cases}\]  It is easy to verify that $G\wr S$ is an
inverse subsemigroup of $\mathfrak I_{G\times [n]}\cong \mathfrak
I_{|G|n}$.  The reader should consult the text of
Eilenberg~\cite{Eilenberg} for more on wreath products of partial
transformation semigroups. We remark that the representation in the
form $(f,\sigma)$ is not unique if $\sigma$ is not totally defined:
only the values of $f$ on $\dom \sigma$ are relevant. Alternatively,
$G\wr S$ can be described as all matrices with entries in $G\cup
\{0\}$ that can be obtained by replacing ones in the rook matrices
from $S$ by arbitrary elements of $G$.  So, for example, the signed
symmetric inverse monoid $\ZZ/2\ZZ\wr \mathfrak I_n$ consists of all
signed rook matrices.

Suppose now that $S\leq \mathfrak I_n$ contains all the idempotents
of $\mathfrak I_n$.  Then the idempotents of $G\wr S$ are precisely
the identities at the subsets of the form $G\times X$ with
$X\subseteq [n]$. Hence, $E(G\wr S)\cong \mathfrak B_n$ and so if
$Y\subseteq X$, then
\[\mu(1_{G\times Y},1_{G\times X}) = (-1)^{|X|-|Y|}.\] The maximal
subgroup of $G\wr S$ at $G\times X$ is isomorphic to the wreath
product $G\wr G_X$ where $G_X$ is the maximal subgroup of $S$ at
$X$.  In particular, for $G\wr \mathfrak I_n$, we see that the
maximal subgroup at $G\times [r]$ is $G\wr \mathfrak S_r$.  So for
the signed symmetric inverse monoid, the maximal subgroups are the
signed symmetric groups of appropriate ranks.  Clearly if
$X\subseteq [n]$, then the set of subsets of the form $G\times Y$
with $Y\subseteq X$ is closed under relative complement in
$G\times X$. Thus to verify that Proposition~\ref{booleanalgebra}
applies to $G\wr S$, we must show that if $(f,s)$ represents an
element $g$ of $G\wr G_X$, then $\mathrm {Fix}(g)$ is of the form
$G\times Y$ with $Y\subseteq X$. But $\mathrm{Fix}(g) = G\times
(1f\inv\cap \mathrm{Fix}(s))$, which is of the required form.

  Let $S(p,r)$ be the Stirling number of the second
kind~\cite{Stanley}. It is given by
\[S(p,r)=\frac{1}{r!}\sum_{k=0}^r (-1)^{r-k}\binom{r}{k}k^p.\]
Then we obtain the following generalization
of~\cite[Example~3.18]{Solomon}, where $S=\mathfrak I_n$ and $G$ is
trivial.

\begin{Thm}
Let $S\leq \mathfrak I_n$ contain all the idempotents and let $G$ be
a finite group.  Let $\theta^p$ be the character of the $p$-tensor
power of the representation of $G\wr \mathfrak I_n\leq \mathfrak
I_{|G|n}$ by rook matrices. Let $X\subseteq [n]$ with $|X|=r$ and
let $G_X$ be the associated maximal subgroup of $S$.  Let $\chi\in
\mathrm{Irr}(G\wr G_X)$. Then
\begin{equation}\label{wreathcase}
(\chi,\theta^p)_S =
\frac{1}{|G|^{r-p}|G_X|}\mathrm{deg}(\chi)r!S(p,r).
\end{equation}
In particular, for $S=\mathfrak I_n$, we obtain
\begin{equation}\label{wreathsolomon}
(\chi,\theta^p)_S = \frac{1}{|G|^{r-p}}\mathrm{deg}(\chi)S(p,r).
\end{equation}
\end{Thm}
\begin{proof}
Since $\rk {1_{G\times Y}} = |G||Y|$, we calculate, using
Proposition~\ref{booleanalgebra}:
\begin{align*}
(\chi,\theta^p)_S &=
\frac{1}{|G|^r|G_X|}\mathrm{deg}(\chi)\sum_{Y\subseteq
X}(-1)^{|X|-|Y|} (|G||Y|)^p\\ & =
\frac{1}{|G|^{r-p}|G_X|}\mathrm{deg}(\chi)\sum_{k=0}^r(-1)^{r-k}\binom{r}{k}k^p\\
&= \frac{1}{|G|^{r-p}|G_X|}\mathrm{deg}(\chi)r!S(p,r).
\end{align*}
In particular, if $G_X=\mathfrak S_r$,~\eqref{wreathcase} reduces
to~\eqref{wreathsolomon}.
\end{proof}

Specializing to the case that $G$ is trivial we obtain:
\begin{Cor}\label{incase}
Let $S\leq \mathfrak I_n$ with $E(S)=E(\mathfrak I_n)$.  Let
$X\subseteq [n]$ with $|X|=r$ and let $G_X$ be the maximal subgroup
of $S$ with identity $1_X$.  Let $\chi$ be an irreducible character
of $G_X$.  Then the multiplicity of $\chi^*$ in the $p$-tensor power
$\theta^p$ of the rook matrix representation of $S$ is given by
\[(\chi,\theta^p)_S= \frac{1}{|G_X|}\mathrm{deg}(\chi)r!S(p,r).\]  In particular, if
$S=\mathfrak I_r$, $X=[r]$ and $\chi^{\lambda}$ is the irreducible
character corresponding to a partition $\lambda$ of $[r]$, then the
multiplicity of $(\chi^{\lambda})^*$ in $\theta^p$ is
$f_{\lambda}S(p,r)$ where $f_{\lambda}$ is the number of standard
Young tableaux of type $\lambda$ (c.f.~\cite{sagan}).
\end{Cor}

 Two other examples of semigroups for which the formula from
Corollary~\ref{incase} is valid are $\mathfrak B_n$ and for the
inverse semigroup  of all order-preserving partial permutations of
$[n]$. In both these cases all the maximal subgroups are trivial. So
the formula then comes down to saying that if $e$ is an idempotent
of rank $r$ in either of these semigroups, then the multiplicity of
the unique irreducible representation associated to $e$ in
$\theta^p$ is simply $r!S(p,r)$.

\subsection{Exterior powers}
We continue to take $K$ to be an algebraically closed field of
characteristic zero.   Solomon showed~\cite[Example~3.22]{Solomon}
that the exterior powers of the rook matrix representation of
$\mathfrak I_n$ are irreducible and are induced by the alternating
representations of the maximal subgroups.  The proof is quite
involved: he makes heavy use of the combinatorics of Ferrer's
diagrams and the theory of symmetric functions.  Here we obtain a
more general result with a completely elementary proof.

Let us denote by $\mathfrak A_X$  the alternating group on a
finite set $X$. If $S\leq \mathfrak I_n$ and $G_e$ is a maximal
subgroup with identity $e$, then define, for $g\in G_e$,
\[\sgn g e=
\begin{cases} 1 & g\in \mathfrak A_{\dom e}\\ -1 &
\text{otherwise.}\end{cases}\]  That is, $\sgn g e$ is the sign of $g$ as
a permutation of $\dom e$.  Then the map
\mbox{$\mathrm{sgn}_e:G_e\to K$} is an irreducible representation
(it could be the trivial representation if $G_e\subseteq \mathfrak
A_{\dom e}$).

\begin{Thm}
Let $S\leq \mathfrak I_n$ and let $\theta^{\wedge p}$ denote the
character of the $p^{th}$-exterior power of the rook matrix
representation of $S$.  Suppose that $S$ contains all the rank $p$
idempotents of $\mathfrak I_n$.   Let $G_e$ be a maximal subgroup of $S$
and $\chi\in \mathrm{Irr}(G_e)$. Then
\[(\chi,\theta^{\wedge p})_S = \begin{cases} 1 & \rk e=p, \chi =
\mathrm{sgn}_e\\ 0 & \text{otherwise.}\end{cases}\]  In particular, under
these hypotheses $\theta^{\wedge p}$ is an irreducible character if
and only if all rank $p$ idempotents of $S$ are isomorphic.  So for
$\mathfrak I_n$, all the exterior powers of the rook matrix representation are irreducible.
\end{Thm}
\begin{proof}
Let $V$ be the module affording the rook matrix representation and
$\{e_1,\ldots, e_n\}$ be the standard basis; so a basis for
$\bigwedge^p V$ is the set \[\{e_{i_1}\wedge \cdots \wedge
e_{i_p}\mid 1\leq i_1<\cdots<i_p\leq n\}.\]  If $s\in S$, then
\[(e_{i_1}\wedge \cdots \wedge e_{i_p})s = e_{i_1s}\wedge\cdots
\wedge e_{i_ps}.\]  This will be a non-zero multiple of
$e_{i_1}\wedge \cdots \wedge e_{i_p}$ if and only if $Xs = X$, where
$X=\{e_{i_1},\ldots,e_{i_p}\}$.  In this case, \[(e_{i_1}\wedge
\cdots \wedge e_{i_p})s = \sgn {s|_X} X (e_{i_1}\wedge \cdots \wedge
e_{i_p}).\]  Thus
\begin{equation}\label{altchar}
\theta^{\wedge p}(s) = \sum_{Y\subseteq \dom s, |Y|=p, Ys=Y} \sgn
{s|_Y} Y.
\end{equation}

We now compute for $\chi\in \mathrm {Irr}(G_e)$:
\begin{equation}\label{eq:ext}
\begin{split}
(\chi,\theta^{\wedge p})_S &= \sum_{f\leq e}
\frac{1}{|G_e|}\sum_{g\in
G_e}\chi(g\inv)\theta^{\wedge p}(fg)\mu(f,e)\\
 &= \frac{1}{|G_e|}\sum_{g\in G_e}\chi(g\inv)\sum_{f\leq e}\theta^{\wedge
 p}(fg)\mu(f,e).
 \end{split}
\end{equation}

Using~\eqref{altchar}, we obtain
\begin{align*}
\sum_{f\leq e}\theta^{\wedge p}(fg)\mu(f,e) &= \sum_{f\leq
e}\left(\sum_{Y\subseteq \dom f, |Y|=p, Yfg=Y} \sgn {(fg)|_Y} Y\right)\mu(f,e)\\
& = \sum_{f\leq e}\left(\sum_{Y\subseteq \dom f, |Y|=p, Yg=Y} \sgn
{g|_Y} Y)\right)
 \mu(f,e) \\
&= \sum _{Y\subseteq \dom e, |Y|=p, Yg=Y} \sgn {g|_Y}
Y\left(\sum_{1_Y\leq f\leq e}\mu(f,e)\right).
\end{align*}
An application of~\eqref{mobiusrestate} then shows that
$\sum_{1_Y\leq f\leq e}\mu(f,e)$ is zero unless \mbox{$1_Y=e$}, in
which case it is one.  So if $\rk e\neq p$, $(\chi,\theta^{\wedge
p})_S=0$.  Otherwise,~\eqref{eq:ext} becomes
\[(\chi,\theta^{\wedge p})_S =\frac{1}{|G_e|}\sum_{g\in
G_e}\chi(g\inv)\sgn g e = (\chi,\mathrm{sgn}_e)_{G_e}.\]  This
establishes the first part of the theorem.  The second part is
immediate from the first.
\end{proof}

\subsection{Direct products}
Let us consider a direct product $S\times T$ of two inverse
semigroups $S,T$.  Then $E(S\times T) = E(S)\times E(T)$.  Moreover,
the natural order on $S\times T$ is the product ordering.  Hence,
by~\eqref{mobiusproduct}, \[\mu_{S\times T}((s,t),(s',t')) =
\mu_S(s,s')\mu_T(t,t').\]  The $\D$-relation on $S\times T$ is the
product of the $\D$-relations on $S$ and $T$ and the maximal
subgroup at an idempotent $(e,f)$ is $G_e\times G_f$.  So if: $K$ is
a commutative ring with unit; $e_1,\ldots, e_s$, respectively
$f_1,\ldots, f_t$, represent the $\D$-classes of idempotents of $S,T$,
respectively; and $n_i,m_j$ are the number of
idempotents in the $\D$-classes of $e_i,f_j$, respectively, then
\begin{align*}
K(S\times T) &\cong \bigoplus_{i,j} M_{n_im_j}(K(G_{e_i}\times
G_{f_j}))\\ & \cong \bigoplus_{i,j} M_{n_im_j}(KG_{e_i}\otimes
KG_{f_j})\\ & \cong \bigoplus_{i,j} \left(M_{n_i}(KG_{e_i})\otimes
M_{m_j}(KG_{f_j})\right)\\ & \cong KS\otimes KT.
\end{align*}
If $K$ is a field, then an irreducible character of $G_{e_i}\times
G_{f_j}$ is of the form $\chi\otimes \eta$ where $\chi\in
\mathrm{Irr}(G_{e_i})$, $\eta\in \mathrm{Irr}(G_{f_j})$ and
$(\chi\otimes \eta)(g,h) = \chi(g)\eta(h)$~\cite{Serre}.  So if
$\theta$ is a character of $S\times T$, then (assuming
characteristic zero) the intertwining number $(\chi\otimes \eta,\theta)_{S\times T}$ is given by
\[\frac{1}{|G_{e_i}||G_{f_j}|}\sum_{e\leq e_i, f\leq
f_j}\sum_{g\in G_{e_i}, h\in G_{f_j}}\chi(g\inv)\eta(h\inv)\theta
(eg,fh)\mu_S(e,e_i)\mu_T(f,f_j).\]

\subsection{Decomposing partial permutation representations}
Let $S\leq \mathfrak{I}_n$ be a partial permutation inverse
semigroup of degree $n$.  We wish to decompose the associated rook
matrix representation into irreducible constituents.  If $S$ is not
transitive, then we can clearly obtain a direct sum decomposition in
terms of the transitive components, so we may as well assume that
$S$ is a transitive partial permutation semigroup. One could obtain
the decomposition below from semigroup theory folklore results: one
argues that in this case the rook matrix representation is induced
from the permutation representation of the unique maximal subgroup $G$ of
the $0$-minimal ideal of $S$ via the Sch\"utzenberger representation
and hence by~\cite{RhodesZalc} and~\cite{Lallement} it decomposes
via the representations induced from those needed to decompose the
associated permutation representation of $G$.  We shall prove this
in a more combinatorial way using our multiplicity formula. Let us
first state the result precisely.

\begin{Thm}\label{decompose}
Let $S\leq \mathfrak{I}_n$ be a transitive partial permutation
inverse semigroup. Let $e\in E(S)$ be an idempotent with associated
maximal subgroup $G_e$. Let $\chi\in \mathrm{Irr}(G_e)$ and let
$\theta$ be the character of the rook matrix representation of $S$
(that is the fixed-point character).  Then
\begin{equation}
(\chi,\theta)_S = \begin{cases} (\chi,\theta|_{G_e})_{G_e} &
\text{if $e$ has minimal non-zero rank} \\ 0 &
\text{else.}\end{cases}
\end{equation}
Moreover, all idempotents of minimal  non-zero rank are isomorphic.
\end{Thm}

This theorem says that decomposing $\theta$ corresponds to
decomposing $\theta|_{G_e}$ where $e$ is an idempotent of minimal
non-zero rank.  Before proving the theorem, we consider some of its
consequences.  The following corollary (well known in the semigroup
representation theory community) is quite useful. See the work of
Zalcstein~\cite{zalc} for an analogue.

\begin{Cor}\label{rlmirred}
Suppose that $S\leq \mathfrak I_n$ is a transitive inverse semigroup
containing a rank $1$ transformation.  Then the rook matrix
representation of $S$ is irreducible and is induced from the trivial
representation of the maximal subgroup corresponding to a rank one
idempotent.
\end{Cor}
\begin{proof}
The maximal subgroup at a rank $1$ idempotent $e$ is trivial and the
restriction of the fixed-point character $\theta$ to $e$ is the
character of the trivial representation of $\{e\}$ since $\theta(e)
= \rk e=1$. Theorem~\ref{decompose} immediately gives the result.
\end{proof}

 The above result then leads to the
following corollary identifying the representation of an inverse
semigroup $S$ induced by the trivial representation of a maximal
subgroup.

\begin{Cor}\label{munnreprlm}
Let $S$ be a finite inverse semigroup and let $D$ be a $\D$-class of
$S$.  Then $S$ acts by partial permutations on the idempotents
$E(D)$ of $D$ via conjugation as follows: \[e^s = \begin{cases}
s\inv es & e\leq \dom s\\ \text{undefined} &
\text{else}\end{cases}\] for $e\in E(D)$ and $s\in S$.  The
associated rook matrix representation of $S$ is the irreducible
representation corresponding to the trivial representation of the
maximal subgroup of $D$.
\end{Cor}
\begin{proof}
We first show that the action is well defined.  Suppose $e\in D$ and
$e^s$ is defined.   Then
\[s\inv ess\inv es = s\inv ees = s\inv es\] as $e\leq ss\inv$.  Thus
$e^s\in E(S)$.  Also
\begin{align*}
\dom {es} &= ess\inv e = e\\
\ran {es} &= s\inv ees = s\inv e s,
\end{align*}
showing that $e^s\in D$.  The reader can verify directly that $e\leq
\dom {st}$ if and only if $e\leq \dom s$ and $e^s\leq \dom t$ and
that in this case $(e^s)^t = e^{st}$.

To complete the proof it suffices by Corollary~\ref{rlmirred} to
show that $S$ acts transitively and that some element of $D$ acts as
a rank one partial transformation. If $e,f\in D$, then there exists
$s\in D$ with $\dom s =e$, $\ran s =f$. Hence
\[e^s=s\inv es = s\inv ss\inv s =s\inv s=\ran s =f\] establishing
transitivity. If $e,f\in D$ and $e\neq f$, then $f\nleq e$
(as idempotents in a $\D$-class of a finite inverse semigroup are
incomparable~\cite{Lawson}).  Hence we may conclude that $e$ acts as
the rank one partial identity that fixes $e$ and is undefined
elsewhere. This completes the proof.
\end{proof}

The above representation is  the right letter mapping
representation of Rhodes~\cite{Arbib,qtheor} in disguise and its
irreducibility follows from the results of~\cite{RhodesZalc}. The
direct sum of the above representations over all $\D$-classes is the
celebrated Munn representation~\cite{Lawson}.

The proof of Theorem~\ref{decompose} requires some preliminary
results that form part of the body of semigroup theory folklore. We
begin with the proof of the final statement.

\begin{Prop}\label{0min}
Let $S\leq \mathfrak{I}_n$ be a transitive inverse semigroup.  Let
$e,f$ be idempotents of minimal non-zero rank.  Then $e\D f$.
\end{Prop}
\begin{proof}
If $e=f$, then  there is nothing to prove.  Assume that $e\neq f$.
Choose $x\in \dom e$, $y\in \dom f$.  By transitivity there exists
$s\in S$ such that $xs=y$.  We claim that $\dom {esf}=e$, $\ran
{esf}= f$ and so $e\D f$.  Indeed, $x\in \dom {esf}\subseteq \dom
e$. Minimality then gives that $\dom {esf} = e$.  The argument that
$\ran {esf} = \ran{f}$ is similar.
\end{proof}

\begin{Lemma}\label{cover}
Let $S\leq \mathfrak{I}_n$ be a transitive inverse semigroup. Then
the domains of the minimal rank non-zero idempotents partition
$[n]$.
\end{Lemma}
\begin{proof}
We first prove disjointness.  Let $i\in \dom e\cap \dom f$ with $e$
and $f$ minimal rank non-zero idempotents.  Then $i\in \dom {ef}$ and
so we have that $0<\rk{ef}\leq \rk e,\rk f$.  Thus $e=ef=f$, as required.  Now we
show that any $i\in [n]$ belongs to the domain of some idempotent
$e$ of minimal non-zero rank.   Let $f$ be any idempotent of minimal
non-zero rank and let $j\in \dom f$. By transitivity there is an
element $s\in S$ with $js=i$. Then $i\in \dom{s\inv fs}$ and $0<\rk
{s\inv fs}\leq \rk f$.  Since $s\inv fs$ is an idempotent, this completes the proof.
\end{proof}

\subsubsection{Proof of Theorem~\ref{decompose}} We shall use the
notation from the statement of Theorem~\ref{decompose}.  Define
$h:S\to \mathbb{Z}$ by setting
\begin{equation}\label{defineh}
 h(s) = \sum_{t\leq s} \theta(t)\mu(t,s).
\end{equation}

 We need a key technical lemma.

\begin{Lemma}\label{tech1}
The function $h:S\to \mathbb{Z}$ is given by the formula: \[h(s) =
\begin{cases} \theta(s) & \text{if $s$ has minimal non-zero rank} \\
0 & \text{otherwise.}\end{cases}\]
\end{Lemma}
\begin{proof}
Let $M$ be the set of elements of $S$ of minimal non-zero rank and
let $M^c$ be the complement of $M$. Recall that $\theta(t) =
|\mathrm{Fix}(t)|$ for $t\in S$.  In particular $\theta(0)=0$ (if
$0\in S$). If $s\in M$, then from $\theta(0)=0$ and~\eqref{defineh}
we may deduce $h(s) = \theta(s)$, as desired. So assume $s\notin M$,
i.e.~$s\in M^c$.

  M\"obius inversion
(Theorem~\ref{inversion}) gives us
\begin{equation}\label{tech1eq}
\theta(s) = \sum_{t\leq s} h(t) = \sum_{t\in M,\ t\leq s}
h(t)+\sum_{t\in M^c,\ t\leq s} h(t)
\end{equation}
By the case already handled, for $t\in M$, \[h(t) =
\theta(t)=|\mathrm{Fix}(t)|.\] Let $t,t'\in M$ with $t,t'\leq s$ and
suppose that $\mathrm{Fix}(t)\cap \mathrm{Fix}(t')\neq \emptyset$.
In particular, $\dom t\cap \dom {t'}\neq \emptyset$.  Since $\dom
t$, $\dom {t'}$ are minimal rank non-zero idempotents,
Lemma~\ref{cover} tells us that $\dom t=\dom {t'}$. It then follows
from the fact that $t,t'$ are both restrictions of $s$ that $t=t'$.
Hence distinct $t\in M$ with $t\leq s$ have disjoint fixed-points
sets implying
\begin{equation}\label{quixote}
\sum_{t\in M,\ t\leq s} h(t) = \sum_{t\in M,\ t\leq s}
|\mathrm{Fix}(t)| = \left|\bigcup _{t\in M,\ t\leq s}
\mathrm{Fix}(t)\right|.
\end{equation}
We claim \[\bigcup_{t\in M,\ t\leq s} \mathrm{Fix}(t) =
\mathrm{Fix}(s).\] Clearly the left hand side is contained in the
right hand side. Conversely, if $x\in \mathrm{Fix}(s)$, then by
Lemma~\ref{cover} there is a (unique) minimal non-zero idempotent
$f$ such that $x\in \dom f$. Hence $x\in \mathrm{Fix} (fs)$, $fs\in
M$ and $fs\leq s$. Thus the sum in~\eqref{quixote} is
$|\mathrm{Fix}(s)| = \theta(s)$. Putting this together
with~\eqref{tech1eq}, we obtain
\begin{equation}\label{quixote2}
0 = \sum_{t\in M^c,\ t\leq s} h(t).
\end{equation}
The formula~\eqref{quixote2} is valid for any element $s\in M^c$.
Therefore, M\"obius inversion in $M^c$ (with the induced ordering)
and~\eqref{quixote2} imply $h(s) = 0$ for any $s\in M^c$, as
desired.
\end{proof}

We may now complete the proof of Theorem~\ref{decompose}.  One
calculates
\begin{align*}
(\chi,\theta)_S &= \sum_{f\leq e} (\chi,\theta_f)_{G_e}\mu(f,e) \\
                &= \frac{1}{|G_e|}\sum_{g\in G_e} \chi(g\inv)
                \sum_{f\leq e} \theta(fg)\mu(f,e)\\
                &=\frac{1}{|G_e|}\sum_{g\in G_e} \chi(g\inv)
                \sum_{t\leq g} \theta(t)\mu(t,g)\\
                & = \frac{1}{|G_e|}\sum_{g\in G_e} \chi(g\inv)h(g)\\
                & = (\chi,h|_{G_e})_{G_e}.
                \end{align*}
But $h|_{G_e}=0$ if $e$ is not of minimal non-zero rank and is
$\theta|_{G_e}$ otherwise by Lemma~\ref{tech1}.  This establishes
Theorem~\ref{decompose}.\qed

\section{The Character Table and Solomon's Approach}\label{compare}
In this section we generalize Solomon's approach to
multiplicities~\cite[Section 3]{Solomon}, based on the character
table, to arbitrary finite inverse semigroups.  We shall see that
this approach leads to greater combinatorial difficulties for
computation.
 If $s$ is an element of a finite semigroup, then $s^{\omega}$
denotes its unique idempotent power.  We set
$s^{\omega+1}=ss^{\omega}$. The element $s^{\omega+1}$ belongs to
the unique maximal subgroup of $\langle s\rangle$.  It is well known
that any character $\chi$ of $S$ has the property that
$\chi(s)=\chi(s^{\omega+1})$~\cite{McAlister,RhodesZalc}.  Indeed,
if $s^n=s^{\omega}$, then $s^{2n}=s^n$.  So if $\rho$ is the
representation affording $\chi$, then $\rho(s)$ satisfies the
polynomial $x^n(x^n-1)$.  It follows that $\rho(s)$ and
$\rho(s)^{n+1}$ have the same trace by considering, say, their
Jordan canonical forms over the algebraic closure.  Thus
$\chi(s)=\chi(s^{n+1})= \chi(s^{\omega+1})$.

Let $S$ be a finite inverse semigroup.  Let $e$ and $f$ be
isomorphic idempotents.  Then $g\in G_e$ and $h\in G_f$ are said to
be \emph{conjugate} if there exists $s\in S$ such that $\dom s = e$,
$\ran s =f$ and $g = shs\inv$.  Note that if $e=f$, this reduces to
the usual notion of conjugacy in $G_e$.  Conjugacy is easily
verified to be an equivalence relation.  Let $\chi$ be a character
of $S$ and suppose that $g$ and $h$ are conjugate, say $g=shs\inv$
with $s$ as above.  Then \[\chi(g) = \chi(shs\inv) = \chi (hs\inv s)
= \chi (h)\] since $s\inv s=\ran s =f$.  Here we used that
$\mathrm{tr}(AB) = \mathrm{tr}(BA)$ for matrices.

Finally define, for $s,t\in S$, $s\sim t$ if $s^{\omega+1}$ is
conjugate to $t^{\omega+1}$.  It is easy to see that this is an
equivalence relation, which we call \emph{character equivalence}.
The discussion above shows that the characters of $S$ are constant
on character equivalence classes.   Let $e_1,\ldots, e_n$ represent
the distinct isomorphism classes of idempotents and set
$G_i=G_{e_i}$.   Then there is a bijection between character
equivalence classes and the (disjoint) union of the conjugacy
classes of the $G_i$. Namely, each character equivalence class
intersects one and only one maximal subgroup $G_i$ and it intersects
that subgroup in a conjugacy class.  We use $\ov s$ for the
character equivalence class of $s\in S$.

One can the define (c.f.~\cite{McAlister,Munnchar,RhodesZalc}) the
\emph{character table} $\mathsf C$ of a finite inverse semigroup $S$
to have rows indexed by the irreducible characters of $S$ over the
complex field $\mathbb{C}$ and columns indexed by the character
equivalence classes. The entry in the row of a character $\chi$ and
the column of a character equivalence class $\ov s$ gives the value
of $\chi$ on $\ov s$.  The table $\mathsf C$ is square since both
sides are in bijection with the (disjoint) union of the conjugacy classes of
the $G_i$.  In order to arrange the table in the most convenient way
possible, let us recall that there is a preorder on the idempotents
of any inverse semigroup $S$ defined as follows~\cite{Lawson}:
$e\preceq f$ if there is an idempotent $e'$ with $e\D e'\leq f$.
This is the same as saying there are elements $s,t\in S$ with $e =
sft$.  For a finite inverse semigroup, $e\preceq f$ and $f\preceq e$
if and only if $e\D f$~\cite{Lawson}. In particular, if
$e_1,\ldots,e_n$ are as above, then $\preceq$ induces a partial
order on $\{e_1,\ldots,e_n\}$. Reordering if necessary, we may
assume that $e_i\preceq e_j$ implies $i\leq j$. It is easy to check
(and well known) that if $\chi$ is an irreducible character coming
from a maximal subgroup $G_i$ and $g\in G_j$ is such that
$\chi(g)\neq 0$, then $e_i\preceq e_j$. Indeed, for the character
not to vanish, $g$ must have a restriction in the $\D$-class of
$e_i$.

Instead of labelling the rows by irreducible characters of $S$ we
label them by $\biguplus \mathrm{Irr}(G_i)$ and, similarly, instead
of labelling the columns by character equivalence classes, we label
them by the conjugacy classes of the $G_i$.  So if $\chi\in
\mathrm{Irr}(G_i)$ and $C$ is a conjugacy class of $G_j$, then
$\mathsf C_{\chi,C}$ is the value of $\chi^*$ on the character
equivalence class containing $C$.  If we group the rows and columns
in blocks corresponding to the ordering $G_1,\ldots,G_n$, then the
character table becomes block upper triangular. More precisely, we
have
\[\mathsf C=
\begin{pmatrix} \mathsf X_1 & \cdots &
* &
*\\ \vdots &\ddots &\vdots &\vdots \\ 0 &\cdots&\mathsf X_{n-1}&*\\
0&\cdots&0&\mathsf X_n\end{pmatrix}\] where $\mathsf X_i$ is the
character table of $G_i$.  In particular, the matrix $\mathsf C$ is
invertible (as character tables of finite groups are invertible).
  Notice
that Solomon's table~\cite{Solomon} differs from ours cosmetically
in how the rows and columns are arranged.

Define a block diagonal matrix
\begin{equation}\label{defineY}
\mathsf Y = \mathrm{diag}(\mathsf X_1,\ldots,\mathsf X_n).
\end{equation}
Then there are unique block upper unitriangular matrices $\mathsf A$
and $\mathsf B$ such that
\begin{equation}\label{defineAB}
\mathsf C = \mathsf Y\mathsf A\qquad \text{and}\qquad \mathsf
C=\mathsf B\mathsf Y.
\end{equation}

So to determine the character table of $S$, one just needs $\mathsf
Y$ (that is the character tables of the maximal subgroups) and
$\mathsf A$ or $\mathsf B$.  We aim to show that $\mathsf A$ is
determined by combinatorial data associated to $S$.  Solomon
explicitly calculated this matrix for $\mathfrak I_n$, but it seems
to be a daunting task in general.  If $g\in G_i$, we use $C^i_g$ to
denote the conjugacy class of $g$ in $G_i$.

\begin{Prop}\label{describeA}
Let $h\in G_i$ and $g\in G_j$.  Then $\mathsf A_{C^i_h,C^j_g}$ is
the number of restrictions of $g$ that are conjugate to $h$ in $S$.
\end{Prop}
\begin{proof}
Let $g\in G_j$ and $\chi$ be an irreducible character coming from
$G_i$.  For each idempotent $f\D e_i$, choose $p_f$ with $\dom {p_f}
= e_i$, $\ran {p_f}=f$.  Given $h\in G_i$, let $a_{C^i_h,C^j_g}$ be
the number of restrictions of $g$ conjugate to $h$ in $S$ (this
number can easily be verified  to depend only on $C^i_h$ and
$C^j_g$).

Then, by the results of~\cite{mobius1} or direct calculation,
\begin{align*}
\mathsf C_{\chi,C^j_g} &= \sum_{f\leq e_j, f \mathrel{\mathscr D}
e_i, g\inv fg=f}\chi(p_f(fg)p_f\inv)\\ & =
\sum_{C^i_h}\chi(h)a_{C^i_h,C^j_g}
\end{align*}
where the last sum is over the conjugacy classes $C^i_h$ of $G_i$.
But
\[\sum_{C^i_h}\chi(h)a_{C^i_h,C^j_g} = \sum_{C^i_h} \mathsf
Y_{\chi,C^i_h}a_{C^i_h,C^j_g}.\]  As $\mathsf Y$ is invertible, a
quick glance at~\eqref{defineAB} allows us to deduce the equality
$\mathsf A_{C^i_h,C^j_g}=a_{C^i_h,C^j_g}$, as required.
\end{proof}
Solomon computed $\mathsf A$ explicitly for $\mathfrak I_n$.  He
showed that if the conjugacy class of $\sigma\in \mathfrak S_r$
corresponds to the partition $\alpha$ of $r$ and the conjugacy class
of $\tau\in \mathfrak S_{\ell}$ corresponds to the partition $\beta$
of $\ell$ where $\alpha$, $\beta$ have respectively $a_i,b_i$ parts
equal to $i$, then
\[\mathsf A_{C_\tau^{\ell},C_\sigma^r} = \binom {\alpha}{\beta} =
\prod_{i\geq 1}\binom{a_i}{b_i}.\]

Similarly we can calculate $\mathsf B$.  If $h\in G_i$, denote by
$z_h$ the size of the centralizer of $h$ in $G_i$;  so the equality
$z_h\cdot |C_h^i|=|G_i|$ holds.

\begin{Prop}\label{describeB}
If $\chi$ is an irreducible character of $G_i$ and $\theta$ is an irreducible character of
$G_j$, then
\begin{equation}\label{bformula}
\mathsf B_{\chi,\theta} = \sum_{C^i_h,C^j_g}z_{g}\inv \mathsf
A_{C^i_h,C^j_g}\chi (C^i_h)\theta(C^j_g)
\end{equation}
 where $C^i_h$ runs over the conjugacy
classes of $G_i$ and $C^j_g$ runs over the conjugacy classes of
$G_j$.
\end{Prop}
\begin{proof}
Define for $1\leq i\leq n$ the matrix $\mathsf Z_i=\mathsf
X_i^{T}\mathsf X_i$.  By the second orthogonality relation for group
characters, $\mathsf Z_i$ is a diagonal matrix whose entry in the
diagonal position corresponding to a conjugacy class $C^i_h$ of
$G_i$ is $z_h$.   Let $\mathsf W = \mathrm{diag}(\mathsf Z_1,\ldots
\mathsf Z_n)$.  Then $\mathsf Y^T\mathsf Y = \mathsf W$. So
from~\eqref{defineAB}, we see that \[\mathsf B = \mathsf Y\mathsf
A\mathsf Y\inv = \mathsf Y\mathsf A\mathsf W\inv \mathsf Y^T.\]
Comparing the $\chi,\theta$ entry of the right hand side of the
above equation with the right hand side of~\eqref{bformula}
completes the proof.
\end{proof}

For $\mathfrak I_n$, Solomon gave a combinatorial interpretation for
the entries of $\mathsf B$ in terms of Ferrer's
diagrams~\cite[Proposition 3.11]{Solomon}. We now turn to the
analogue of Solomon's multiplicity formulas for $\mathfrak
I_n$~\cite[Lemma 3.17]{Solomon}, in terms of $\mathsf A$ and
$\mathsf B$, for the general case. It is not clear how usable these
formulas are since one needs quite detailed information about the
inverse semigroup $S$ to determine these matrices.  Our previous
computations with the multiplicity
formula~\eqref{semigroupintertwine} often just used knowledge of the
idempotent set, while Solomon's approach requires much more. It also
explains why partitions, Ferrer's diagrams and symmetric functions
come into Solomon's approach for the tensor and exterior powers of
the rook matrix representation of $\mathfrak I_n$, but they play no
role in our approach.  We retain the above notation.

\begin{Thm}
Let $\chi$ be an irreducible character of $G_i$ and $\theta$ a
character of $S$.  Then the following two formulas are valid:
\begin{equation}\label{solomonformula1}
(\chi,\theta)_S = \sum_{C^i_h}\chi(h)z_h\inv\sum_{C^j_g}\mathsf
A\inv_{C^j_g,C^i_h}\theta (g).
\end{equation}
where the first sum is over the conjugacy classes $C^i_h$ of $G_i$
and the second is over all conjugacy classes $C^j_g$ of all the
$G_j$; and
\begin{equation}\label{solomonformula2}
 (\chi,\theta)_S = \sum_{j=1}^n
\sum_{\psi\in \mathrm{Irr}(G_j)}\mathsf B_{\psi,\chi}\inv
(\psi,\theta|_{G_j})_{G_j}
\end{equation}
\end{Thm}
\begin{proof}
For an irreducible character $\psi$ of $G_j$, denote by $m_{\psi}$
the multiplicity of $\psi^*$ in $\theta$.  So $\theta = \sum_{\psi}
m_{\psi}\psi^*$.  Hence, for a conjugacy class $C^j_g$ of $G_j$,
\[\theta (g) = \sum_{\psi} m_{\psi}\psi^*(g)= \sum_{\psi}
m_{\psi}\mathsf C_{\psi,C^j_g}\] where the sum runs over all the
irreducible characters of all the $G_j$.  Using this, we obtain:
\begin{equation}\label{solomonhelper}
(\chi,\theta)_S = \sum_{\psi} m_{\psi}\delta_{\psi,\chi} =
\sum_{\psi} m_{\psi}\sum_{C^j_g} \mathsf C_{\psi,C^j_g}\mathsf
C_{C^j_g,\chi}\inv = \sum_{C^j_g}\mathsf C_{C^j_g,\chi}\inv\theta(g)
\end{equation}
where the last sum runs over the conjugacy classes $C^j_g$ of the
$G_j$, \mbox{$j=1,\ldots,n$}. Setting $\mathsf W= \mathsf Y^T\mathsf
Y$ again (as in the proof of Proposition~\ref{describeB}), we
compute:
\[\mathsf C_{C^j_g,\chi}\inv = (\mathsf A\inv \mathsf W\inv \mathsf
Y^T)_{C^j_g,\chi} = \sum_{C^i_h} \mathsf
A\inv_{C^j_g,C^i_h}z_h\inv\chi(h)\] where the last sum runs over the
conjugacy classes $C^i_h$ in $G_i$.  This, in conjunction
with~\eqref{solomonhelper} implies~\eqref{solomonformula1}.  Also,
for $g\in G_j$,
\begin{equation}\label{hamburgerhelper}
\mathsf C_{C^j_g,\chi}\inv = (\mathsf W\inv \mathsf Y^T\mathsf
B\inv)_{C^j_g,\chi} = z_g\inv \sum_{\psi\in \mathrm{Irr}(G_j)}
\psi(g)\mathsf B_{\psi,\chi}\inv.
\end{equation}
  This last equality uses that
$\mathsf Y^T$ is block diagonal.   Combining~\eqref{hamburgerhelper}
with~\eqref{solomonhelper} and the fact that as $C^j_g$ runs over
all conjugacy classes of $G_j$,
\[\sum_{C^j_g}z_g\inv \psi(g)\theta(g) = (\ov
{\psi},\theta)_{G_j},\] where $\ov {\psi}$ is the conjugate
character of $\psi$, gives~\eqref{solomonformula2}.  This completes
the proof of the theorem.
\end{proof}

\section{Semigroups with Commuting Idempotents and Generalizations}
There is a large class of finite semigroups whose representation
theory is controlled in some sense by an inverse semigroup.  For
instance, every irreducible representation of an idempotent
semigroup factors through a semilattice; for a readable account of
this classical fact for non-specialists, see~\cite{Brown2}. This
section will require a bit more of a semigroup theoretic background
than the previous ones; the reader is referred to~\cite{CP} for
basics about semigroups and~\cite{Arbib,qtheor} for results specific to
finite semigroups.  The book of Almeida~\cite{Almeidabook} contains
more modern results, as does the forthcoming book~\cite{qtheor}.

\subsection{Semigroups with commuting idempotents}
Let us first begin with a class that is very related to inverse
semigroups: semigroups with commuting idempotents.  Every inverse
semigroup has commuting idempotents, as does every subsemigroup of
an inverse semigroup.  However, not every finite semigroup with
commuting idempotents is a subsemigroup of an inverse semigroup. It
is a very deep result of Ash~\cite{AshInv,ash} that every finite
semigroup with commuting idempotents is a quotient of a subsemigroup
of a finite inverse semigroup (Ash's original proof uses Ramsey
theory in an extremely clever way~\cite{AshInv}; this result can
also be proved using a theorem of Ribes and Zalesski\u\i\ about the
profinite topology on a free group~\cite{RZ,HMPR}).

Elements $s,t$ of a semigroup $S$ are said to be \emph{inverses} if
$sts=s$ and $tst=t$.  Elements with an inverse are said to be (von
Neumann) \emph{regular}.  Denote by $R(S)$ the set of regular
elements of $S$.   A semigroup in which all elements are regular is
called, not surprisingly, \emph{regular}.  Regular semigroups are
very important: a connected algebraic monoid with zero has a
reductive group of units if and only if it is
regular~\cite{Putcha,Rennerbook}; the semigroup algebra of a finite
regular semigroup is quasi-hereditary~\cite{quiver}.   It is known
that a semigroup $S$ is an inverse semigroup if and only if it is
regular and has commuting idempotents~\cite{Lawson,CP}. More
generally it is known that if the idempotents of a semigroup
commute, then the regular elements have unique inverses~\cite{CP}.
We give a proof for completeness.

\begin{Prop}\label{uniqueinverse}
Suppose $S$ has commuting idempotents and $t,t'$ are inverses of
$s$.  Then $t=t'$.
\end{Prop}
\begin{proof}
Using that $st,ts,st',t's$ are idempotents we obtain:
\begin{align*}
t &= tst = ts(t'st')st = t'stst'st\\   &= t's(tst)st' = t'st' =t'.
\end{align*}
This establishes the uniqueness of the inverse.
\end{proof}

If $S$ is a semigroup with commuting idempotents and $u\in R(S)$,
then we denote by $u\inv$ the (unique) inverse of $u$.  The
following is a standard fact about semigroups with commuting
idempotents.

\begin{Prop}\label{inversesub}
Let $S$ be a finite semigroup with commuting idempotents.  Then the
set $R(S)$ of regular elements of $S$ is an inverse semigroup.
\end{Prop}
\begin{proof}
The key point is that $R(S)$ is a subsemigroup of $S$.  Indeed, if
\mbox{$a,b\in R(S)$}, then we claim that $(ab)\inv = b\inv a\inv$.
Observing that $aa\inv$, $a\inv a$, $bb\inv$ and $b\inv b$ are
idempotents, we compute: \[ab(b\inv a\inv)ab = a(bb\inv)(a\inv a)b =
a(a\inv a)(bb\inv)b =ab.\] Similarly one verifies $b\inv a\inv
(ab)b\inv a\inv = b\inv a\inv$. Since the inverse of any regular
element is regular, $R(S)$ is closed under taking inverses and hence
is a regular semigroup with commuting idempotents and therefore an
inverse semigroup.
\end{proof}

If $A$ is a finite dimensional algebra, then $\mathrm{Rad}(A)$ denotes the Jacobson radical of $A$.
Our next goal is to show that if $S$ has commuting idempotents and
$K$ is a field of characteristic zero, then $KS/\mathrm{Rad}(KS)
\cong KR(S)$ and the isomorphism is the identity on $KR(S)$ (viewed as a subalgebra of $KS$).  This
means that we can use our results for inverse semigroups to obtain
character formulas for multiplicities of irreducible constituents in
representations of $S$.

\begin{Prop}\label{orderec}
Let $S$ be a semigroup with commuting idempotents. Let $u\in R(S)$
and $s\in S$. Then the following are equivalent:
\begin{enumerate}
\item $uu\inv s = u$;
\item $u= es$ with $e\in E(S)$;
\item $su\inv u =u$;
\item $u = sf$ with $f\in E(S).$
\end{enumerate}
\end{Prop}
\begin{proof}
Clearly (1) implies (2).  For (2) implies (1), suppose that $u=es$
with $e\in E(S)$.  Then $eu=ees=es=u,$ so \[u=uu\inv u = uu\inv es =
euu\inv s = uu\inv s.\]  The equivalence of (3) and (4) is dual.

To prove that (1) implies (3), assume $uu\inv s = u$.  We show that
$su\inv u$ is an inverse for $u\inv$.  Then the equality $su\inv
u=u$ will follow by uniqueness of inverses in $R(S)$
(Proposition~\ref{inversesub}). Indeed
\[(su\inv u)u\inv (su\inv u) = su\inv (uu\inv s) u\inv u = su\inv u
u\inv u = su\inv u.\]  Also, using $u\inv = u\inv uu\inv$,  \[u\inv
(su\inv u)u\inv = u\inv (uu\inv s)u\inv = u\inv uu\inv  = u\inv.\]
The implication (3) implies (1) is proved similarly.
\end{proof}

Let $S$ be a semigroup with commuting idempotents.  Define, for
$s\in S$, \[s^{\downarrow} = \{u\in R(S)\mid uu\inv  s = u\}.\]  In
other words, $s^{\downarrow}$ is the set of elements for which the
equivalent conditions of Proposition~\ref{orderec} hold.   Notice
that if $s$ is regular, then $s^{\downarrow}$ just consists of all
elements of $R(S)$ below $s$ in the natural partial order on $R(S)$, whence
the notation. We can now establish an analogue of
Lemma~\ref{groupoid=semigroup}.

\begin{Lemma}\label{mapin}
Let $\nu:S\to KR(S)$ be given by $\nu(s) = \sum_{t\in
s^{\downarrow}} \lfloor t\rfloor$.  Then $\nu$ is a homomorphism that restricts to the identity on $R(S)$, where we view $R(S)$ as a subsemigroup of $KR(S)$.
\end{Lemma}
\begin{proof}
The observation before the proof shows that if $s\in R(S)$, then
\[\nu(s) = \sum_{t\leq s} \lfloor t\rfloor = s\] via M\"obius
inversion and Theorem~\ref{mainthm:mobius}.

For arbitrary $s,t\in S$, we have that \[\nu(s)\nu(t) = \sum_{u\in
s^{\downarrow}, v\in t^{\downarrow},\ \ran u=\dom v} \lfloor
uv\rfloor.\]  First we show that if $u\in s^{\downarrow}$, $v\in
t^{\downarrow}$ and $\ran u=\dom v$ then $uv\in (st)^{\downarrow}$.
Indeed $\dom {uv}=\dom u$.  Hence $uv(uv)\inv = uu\inv$.  Also
$u\inv u=vv\inv$ so
\[(uv)(uv)\inv st = uu\inv st = ut = u(u\inv u)t = uvv\inv t = uv.\]  To
complete the proof, we must show that every element $u$ of
$(st)^{\downarrow}$ can be written uniquely in the form $s't'$ with
$s'\in s^{\downarrow}$, $t'\in t^{\downarrow}$ and $\ran {s'}=\dom
{t'}$. The proof proceeds exactly along the lines of that of
Lemma~\ref{groupoid=semigroup}. Namely, to obtain such a
factorization $u=s't'$ we must have $\dom{s'}=\dom u$ and $\ran
{t'}=\ran u$, that is, we are forced to set $s'=uu\inv s$ and
$t'=tu\inv u$ since, by Proposition~\ref{orderec}, $s' = s'(s')\inv
s$ and $t'=t(t')\inv t'$. Then one shows, just as in the proof of
Lemma~\ref{groupoid=semigroup}, that $s'$ is regular with inverse
$tu\inv$ and $t'$ is regular with inverse $u\inv s$. Hence $s'\in
s^{\downarrow}$ and $t'\in t^{\downarrow}$ (the latter requires
Proposition~\ref{orderec}). The proofs that $s't'=u$ and $\ran
{s'}=\dom {t'}$ also proceed along the lines of the proof of
Lemma~\ref{groupoid=semigroup} and so we omit them.
\end{proof}

Our next task is to show that the induced surjective homomorphism \mbox{$\nu:KS\to KR(S)$} has
nilpotent kernel. In any event $\nu$ splits as a $K$-vector space
map and so $KS = KR(S)\oplus \ker \nu$ as $K$-vector spaces. Clearly
then $\ker \nu$ has basis $B= \{s-\nu(s)\mid s\in S\setminus
R(S)\}$. Indeed, the number of elements in $B$ is the dimension of
$\ker \nu$ and these elements are clearly linearly independent since
the unique non-regular element in the support of $s-\nu(s)$ is $s$
itself. Notice that, via $\nu$, any irreducible representation of
$R(S)$ extends to $S$. We shall prove the converse.  Our method of
proof will show that $\ker \nu$ is contained in the radical of $KS$.
First we need some definitions.

An \emph{ideal} of a semigroup $S$ is a subset $I$ such that
$SI\cup IS\subseteq I$.   We then place a preorder on $S$ by ordering
elements in terms of the principal ideal they generate.  That is, if
$s\in S$, let $J(s)$ be the principal ideal generated by $s$. Define
$s\leq_{\J} t$ if $J(s)\subseteq J(t)$.  We write $s\J t$ if
$J(s)=J(t)$.  This is an example of one of Green's
relations~\cite{Green,CP,Arbib,qtheor}.  There are similar relations,
denoted $\R$ and $\L$, corresponding to principal right and left
ideals, respectively.  The $\J$-relation on a finite inverse
semigroup $S$ coincides with what we called $\D$
earlier~\cite{Lawson}.  That is, for $s,t\in S$, we have $s\J t$ if
and only if $\dom s$ is isomorphic to $\dom t$. It is known that in
a finite semigroup the following are equivalent for a $\J$-class
$J$~\cite{CP,Arbib,qtheor}:
\begin{itemize}
\item $J$ contains an idempotent;
\item $J$ contains a regular element;
\item Every element of $J$ is regular.
\end{itemize}
Such a $\J$-class is called a \emph{regular $\J$-class}.   Any
inverse of a regular element belongs to its $\J$-class.   If $s,t$
are $\J$-equivalent regular elements, then there are (regular)
elements $x,y,u,v$ ($\J$-equivalent to $s$ and $t$) such that
$xsy=t$ and $utv =s$~\cite{CP,Arbib,qtheor}.  Hence if $S$ has commuting
idempotents, then the $\J$-classes of $R(S)$ are precisely the
regular $\J$-classes of $S$.  The following proof uses a result of
Munn~\cite{Munn} that is exposited in~\cite{CP}.
See~\cite{RhodesZalc,rhodeschar,myrad} for refinements.

\begin{Thm}\label{radialEC}
Let $S$ be a finite semigroup with commuting idempotents and let
$\rho:S\to M_n(K)$ be an irreducible representation. Then:
\begin{enumerate}
\item $\rho|_{R(S)}$ is an irreducible representation;
\item $\ker \nu\subseteq \ker \rho$.
\end{enumerate}
\end{Thm}
\begin{proof}
Let $J$ be a $\leq_{\J}$-minimal $\J$-class of $S$ on which $\rho$
does not vanish.  It is a result of Munn~\cite[Theorem 5.33]{CP}
(see also~\cite{myrad,RhodesZalc}) that $J$ must be a regular
$\J$-class, say it is the $\J$-class of an idempotent $e$ and that
$\rho$ must vanish on any $\J$-class that is not $\leq_{\J}$-above
$J$~\cite{CP,RhodesZalc,myrad}.   Let $I=J(e)\setminus J$. Then $I$
is an ideal of $J(e)$ on which $\rho$ vanishes and so $\rho|_{J(e)}$
factors through a representation \ov {\rho} of the quotient
$J^0=J(e)/I$. Munn proved~\cite[Theorem 5.33]{CP} that $\ov {\rho}$
is an irreducible representation of $J^0$. Since $J$ is a regular
$\J$-class and $J^0=J(e)/I= (J(e)\cap R(S))/(I\cap R(S))$, it
follows that $\ov {\rho}$ is also induced by $\rho|_{J(e)\cap
R(S)}$. It is then immediate that $\rho|_{R(S)}$ is an irreducible
representation since any $R(S)$-invariant subspace is $(J(e)\cap
R(S))$-invariant. Moreover, since $\rho|_{R(S)}$ is induced by $\ov
{\rho}$, it must be an irreducible representation of $KR(S)$
associated to the direct summand of $KR(S)$ spanned by $\{\lfloor
s\rfloor\mid s\in J\}$, which we denote $KJ$ (recall that $J$ is a
$\D$-class).

Now let $s\in S$.  We show that $\rho(s)=\rho(\nu(s))$.  Since
elements of the form $s-\nu(s)$ span $\ker \nu$, this will show that
$\ker \nu\subseteq \ker \rho$.  If $s$ is not $\leq_{\J}$-above $J$,
then $\rho(s)$ is zero.  Since each element of $s^{\downarrow}$ is
$\leq_{\J}$-below $s$, we also have in this case that $\rho(\nu(s))=0$.
Suppose now that $s\geq_{\J} J$, i.e.\ $J\subseteq J(s)$.   Consider $1_J=\sum_{f\in E(J)} \lfloor f\rfloor$. This is the
identity element of $KJ$ and hence is sent to the identity matrix
under $\rho$, as $\rho|_{R(S)}$ is induced by first projecting to
$KJ$.  Therefore, $\rho(s) = \rho (1_Js)$.  It thus suffices to show
that $\rho (1_Js)=\rho (\nu (s))$.  Since, for $t\in J$, every
summand but $t$ of $\lfloor t\rfloor$ is strictly $\leq_{\J}$-below
$J$, we see that in this case $\rho (\lfloor t\rfloor) = \rho (t)$.
Now for $f\in E(J)$ either: $fs<_{\J} J$, and hence $\rho (\lfloor
f\rfloor s)=0$; or $fs\in J$, and so $fs\in s^{\downarrow}$.
Conversely, if $t\in J\cap s^{\downarrow}$, then $t= tt\inv s$ and
$tt\inv\in E(J)$. Thus
\begin{equation}\label{painfuleq1}
\rho (1_Js) = \sum_{f\in E(J), fs\in J} \rho(\lfloor f\rfloor s) =
\sum_{f\in E(J), fs\in J} \rho (fs) = \sum_{t\in s^{\downarrow}\cap
J} \rho (t).
\end{equation}

 Suppose $t\in
s^{\downarrow}$. If $t\notin J$, then $\lfloor t\rfloor$ is not in
$KJ$ and so $\rho(\lfloor t\rfloor) =0$.  Thus
\begin{equation}\label{painfuleq2}
\rho (\nu(s)) = \sum_{t\in s^{\downarrow}\cap J} \rho(\lfloor
t\rfloor) = \sum_{t\in s^{\downarrow}\cap J} \rho (t).
\end{equation}
Comparing~\eqref{painfuleq1} and~\eqref{painfuleq2} shows that
$\rho(1_Js) = \rho (\nu (s))$, establishing that $\ker \nu\subseteq
\ker \rho$.
\end{proof}

\begin{Cor}\label{commutingcase}
Let $S$ be a finite semigroup with commuting idempotents and $K$ a
field.  Define $\nu:KS\to KR(S)$ on $s\in S$ by \[\nu(s) =
\sum_{t\in s^{\downarrow}}\lfloor t\rfloor.\]  Then $\nu$ is a
retraction with nilpotent kernel.  Hence we have the equality
$KS/\mathrm{Rad(KS)}=KR(S)/\mathrm{Rad}(KR(S))$.

In particular, if the characteristic of $K$ is $0$ (or more
generally if the characteristic of $K$ does not divide the order of
any maximal subgroup of $S$), then $\ker \nu=\mathrm{Rad}(KS)$. In
this case $\mathrm{dim}(\mathrm{Rad}(KS)) = |S\setminus R(S)|$ and a
basis for $\mathrm{Rad}(KS)$ is given by the set $\{s-\nu(s)\mid
s\in S\setminus R(S)\}.$
\end{Cor}
\begin{proof}
Theorem~\ref{radialEC} shows that $\ker \nu$ is contained in the
kernel of every irreducible representation of $KS$ and hence $\ker
\nu$ is a nilpotent ideal.  From this the first paragraph follows.
In the context of the second paragraph, we have that $KR(S)$ is
semisimple and so has no nilpotent ideals.  Thus $\ker \nu$ is the
largest nilpotent ideal of $KS$ and hence is the radical.  The
remaining statements are clear.
\end{proof}

It follows directly that the irreducible representations of $S$ are in
bijection with the irreducible representations of its maximal
subgroups up to $\J$-equivalence (actually this is true for any
finite semigroup~\cite{CP,RhodesZalc}).  Moreover, our multiplicity
formulas for inverse semigroups apply verbatim for semigroups $S$
with commuting idempotents.

\begin{Thm}\label{multforEC}
Let $S$ be a finite semigroup with commuting idempotents and $K$ a
field of characteristic zero.   Let $\chi$ be an irreducible
character of a maximal subgroup $G$ with identity $e$ and let $d$ be
the dimension of the associated endomorphism division algebra.  Then
if $\theta$ is a character of $S$ and $m$ is the multiplicity of the
irreducible representation of $S$ associated to $\chi$ as a
constituent in $\theta$, then
\[md=\sum_{f\leq e}(\chi,\theta_f)_G\mu(f,e)\] where $\mu$ is the
M\"obius function of $E(S)$.
\end{Thm}

It follows directly that a completely reducible representation of a finite
semigroup with commuting idempotents is determined by its character.
In general the irreducible constituents are determined by the
character. (Of course this is true in any finite dimensional algebra.)

\subsection{Multiplicities for more general classes of semigroups}
We now want to consider a wider class of semigroups whose
irreducible representations are controlled by inverse semigroups.
First we recall the notion of an \pv{LI}-morphism, which is the
semigroup analogue of an algebra homomorphism with nilpotent kernel.
A finite semigroup $S$ is said to be \emph{locally trivial} if, for
each idempotent $e\in S$, $eSe=e$. A homomorphism $\p:S\to T$ is
said to be an \emph{$\pv{LI}$-morphism} if, for each locally trivial
subsemigroup $U$ of $T$, the semigroup $\p\inv (U)$ is again locally
trivial. The following result, showing that $\pv{LI}$-morphisms
correspond to algebra morphisms with nilpotent kernel, was proved
in~\cite{myrad}.

\begin{Thm}\label{radicalthm}
Let $K$ be a field and $\p:S\to T$ a homomorphism of finite
semigroups.  If $\p$ is an $\pv{LI}$-morphism, then the induced map
\mbox{$\ov{\p}:KS\to KT$} has nilpotent kernel.  The converse holds
if the characteristic of $K$ is zero.
\end{Thm}

Therefore, if $S$ is a finite semigroup with an $\pv{LI}$-morphism
$\p:S\to T$ to a semigroup $T$ with commuting idempotents, then we
can conclude that $KS/\mathrm{Rad}(KS) = KR(T)/\mathrm{Rad}(KR(T))$
(equals $KR(T)$ if $\mathrm{char}(K)=0$).  In particular we can use
our multiplicity formula for inverse semigroups to calculate
multiplicities for irreducible constituents for representations of
$S$ in characteristic zero.  For instance, if $S$ is an idempotent
semigroup, then one can find such a map $\p$ with $T$ a semilattice.
This is what underlies part of the work of
Brown~\cite{Brown1,Brown2}, as well as some more general work of
Putcha~\cite{quiver}. See also~\cite{mobius1}.

Let us describe those semigroups with such a map $\p$.  This class
is well known to semigroup theorists and it would go too far afield
to give a complete proof here, so we restrict ourselves to just
describing the members of the class.  First we describe the
semigroups with an $\pv {LI}$-morphism to a semilattice $L$.  This
class was first introduced by Sch\"utzenberger~\cite{Schutzambi} in
the context of formal language theory.  It consists precisely of
those finite semigroups $S$ such that $R(S)=E(S)$, that is those
finite semigroups all of whose regular elements are idempotents.
This includes of course all idempotent semigroups.  The semilattice
$L$ is in fact the set $\mathscr U(J)$ of regular $\J$-classes
ordered by $\leq_{\J}$. The map sends $s\in S$ to the $\J$-class of
its unique idempotent power. See~\cite{mobius1} for details.  The
class of such semigroups is usually denoted \pv{DA} in the semigroup
literature (meaning that regular $\D$-classes are aperiodic
subsemigroups). It was shown in~\cite{myrad} that \pv{DA}
consists precisely of those finite semigroups with a faithful upper triangular matrix representation over a field of characteristic $0$ by matrices whose only eigenvalues are~$0$ and~$1$.

Now if $\p:S\to T$ is an $\pv {LI}$-morphism, then it is not to hard
to show that if $U\leq S$ is a subsemigroup, then $\p|_U$ is again
an $\pv{LI}$-morphism.  Suppose that $\p:S\to T$ is an $\pv
{LI}$-morphism of finite semigroups where $T$ is a semigroup with
commuting idempotents.  Then $E(T)$ is a semilattice and $\langle
E(S)\rangle$ maps onto $E(T)$ via $\p$.  Hence $\langle E(S)\rangle
\in \pv {DA}$.  Let us denote by \pv{EDA} the collection of all
finite semigroups $S$ such that $E(S)$ generates a semigroup in
\pv{DA}; this includes all semigroups whose idempotents form a
subsemigroup (in particular the class of so-called orthodox
semigroups~\cite{CP}). It is well known to semigroup theorists that
if $S\in \pv {EDA}$, then $S$ admits an \pv{LI}-morphism $\p:S\to T$
to a semigroup $T$ of partial permutations.  In particular, $T$ has
commuting idempotents.  The transitive components of $T$ correspond
to the action of $S$ on the right of a regular $\R$-class of $S$
modulo a certain equivalence relation corresponding to identifying
elements that differ by right multiplication by an element of the
idempotent-generated subsemigroup; the reader can look
at~\cite{Therien} to infer details. Alternatively, one can easily
verify that each generalized group mapping image of $S$
corresponding to a regular $J$-class acts by partial permutations on
its $0$-minimal ideal~\cite{Arbib,qtheor}.   The key point is that $\p$ is
explicitly constructible and hence the multiplicity formulas for
calculating irreducible constituents for representations of $T$ can
be transported back to $S$.  The congruence giving rise to $\p$ is
defined as follows.  Let $S\in \pv{EDA}$ and $s,t\in S$.  Define
$s\equiv t$ if, for each regular $\J$-class $J$ of $S$ and each
$x,y\in J$, one has either $xsy=xty$ or both $xsy,xty\notin
J$~\cite{Arbib,myrad,rhodeschar,qtheor}.
 The details are left to
the reader.

Just to give a sample computation, let $S\in \pv {DA}$ and suppose
that the map \mbox{$\p:S\to M_n(K)$} is an irreducible
representation. The maximal subgroups of $S$ are trivial.  If $J$ is
a regular $\J$-class, then the (unique) irreducible representation
of $S$ associated to $J$ is given by
\[\rho_J(s) =
\begin{cases} 1 & s\geq_{\J} J \\ 0 & \text{else}\end{cases}\] (see~\cite{mobius1}). To
obtain a formula for the multiplicity of $\rho _J$ in $\p$, we must
choose an idempotent $e_J$ for each regular $\J$-class $J$.  Then
the multiplicity of $\rho_J$ in $\p$ is given by
\begin{equation}\label{DAvers}
\sum_{J'\leq_{\J} J, J'\in \mathscr U(J)} \rk
{\p(e_{J'}e_Je_{J'})}\mu(J',J)
\end{equation}
 where $\mu$ is the M\"obius function
of the semilattice $\mathscr U(J)$. This generalizes the
multiplicity results in~\cite{Brown1,Brown2,mobius1} for random
walks on minimal left ideals of semigroups in \pv {DA}.

This yields a direct proof that a completely reducible representation of a semigroup
from $\pv{EDA}$ is determined by its character and that in general
the irreducible constituents are determined by the character.

\subsection{Random walks on triangularizable finite semigroups}
We can now answer a question that remained unsettled
in~\cite{mobius1}. In that paper we calculated the eigenvalues for
random walks on minimal left ideals of finite semigroups admitting a
faithful representation by upper triangular matrices over
$\mathbb{C}$.  This generalized the work of Bidigare \textit{et
al.}~\cite{BHR} and Brown~\cite{Brown1,Brown2}. We showed that there
was an eigenvalue corresponding to each irreducible character of the
semigroup (and gave a formula for the eigenvalue) but at the time we
could only prove that the multiplicity of the eigenvalue was the
same as the multiplicity of the corresponding irreducible
representation as a constituent in the linear representation induced
by the left action on the minimal left ideal. We could only
calculate the multiplicities explicitly if the semigroup belonged to
$\pv {DA}$. With our new tools we can now handle the general case.

So let us call a finite semigroup $S$ \emph{triangularizable} if it
can be represented faithfully by upper triangular matrices over
$\mathbb{C}$.  These were characterized in~\cite{myrad} as precisely
those semigroups admitting an $\pv {LI}$-morphism to a commutative
inverse semigroup. Equivalently, they were shown to be those finite
semigroups in which all maximal subgroups are abelian, whose idempotents generate a subsemigroup with only trivial subgroups and such that
each regular element satisfies an identity of the form $x^m=x$.
Moreover, it was shown that every complex irreducible representation
of such a semigroup has degree one~\cite{myrad} (so its semigroup algebra is basic).  See
also~\cite{mobius1} for more.  Important examples include abelian
groups and idempotent semigroups, including the face semigroup of a
hyperplane arrangement~\cite{Aguiar,Brown1,Brown2}.

Let $S$ be a fixed finite triangularizable semigroup.    Let
$\p:S\to T$ be its \pv {LI}-morphism to a commutative inverse
semigroup.   The semigroup $T$ has a unique idempotent in each
$\D$-class and the corresponding maximal subgroup is abelian.  The
lattice of idempotents of $T$ is isomorphic to the poset $\mathscr
U(S)$ of regular $\J$-classes of $S$.  Fix an idempotent $e_J$ for
each $J\in \mathscr U(S)$.  The maximal subgroup $G_{e_J}$ will be
denoted $G_J$.  We recall the description of the irreducible
characters of $S$. Suppose $G_J$ is a maximal subgroup and $\chi$ is
an irreducible character of $G_J$. Then the associated irreducible
character $\chi^*:S\to \mathbb{C}$ is given by
\[\chi^*(s) = \begin{cases} \chi(e_Jse_J) & s\geq_{\J} J\\ 0
&\text{else.}\end{cases}\] (c.f.~\cite{mobius1}).

Suppose one puts a probability measure $\pi$ on $S$.  That is we
assign probabilities $p_s$ to each $s\in S$ such that $\sum_{s\in
S}p_s=1$. We view $\pi$ as the element $\pi=\sum_{s\in S}p_ss$ of
$\mathbb CS$. Without loss of generality we may assume that $S$ has
an identity; indeed, if $S$ does not have an identity, we can always
adjoin an identity $1$ and set $p_1=0$ without changing the Markov
chain. This assumption guarantees that any representation of $S$
that sends the identity to the identity matrix does not have a null
constituent. Let $L$ be a minimal left ideal of $S$ (what follows is
independent of the choice of $L$ since all minimal left ideals of a
finite semigroup are isomorphic via right translation by Green's
lemma~\cite{CP,Green}). We remark that, for random walks on finite
semigroups, one considers minimal left ideals because if the support
of the probability measure generates the semigroup, then the walk
almost surely enters a minimal left ideal~\cite{Rosenblatt}. The
associated random walk on $L$ is then the Markov chain with
transition operator the $|L|\times |L|$-matrix $M$ that has in entry
$\ell_1,\ell_2$ the probability that if one chooses $s\in S$
according to the probability distribution $\pi$, then
$s\ell_1=\ell_2$.

It is easy to see that if one takes the representation $\rho$
afforded by $\mathbb CL$ (viewed as a left $\mathbb CS$-module),
then $M$ is the transpose of the matrix of
$\rho(\pi)$~\cite{Brown1,Brown2}.  Now since all irreducible
representations of a triangularizable semigroup have degree one, a
composition series for $\mathbb CL$ puts $\rho$ in upper triangular
form with the characters of $S$ on the diagonal, appearing with
multiplicities according to their multiplicities as constituents of
$\rho$. Hence there is an eigenvalue $\lambda_{\chi}$ of $M$
associated to each irreducible character $\chi$ of a maximal
subgroup $G_J$ of $S$ (where $J$ runs over $\mathscr U(S)$), given
by the character sum
\begin{equation}\label{eigenvalue}
\lambda_{\chi} = \chi^*(\pi) =  \sum_{s\in S}
p_s\chi^*(s)=\sum_{s\geq_{\mathscr J} J} p_s\chi (e_Jse_J).
\end{equation}
Of course, it could happen that different
characters yield the same eigenvalue.

  For $s\in S$, let
$\mathrm{Fix}_L(s)$ denote the set of fixed-points of $s$ acting on
the left of $L$. Then the character $\chi_{\rho}$ of $\rho$ simply
counts the cardinality of $\mathrm{Fix}_L(s)$. It is now a
straightforward exercise in applying~\eqref{semigroupintertwine} to
verify that the multiplicity of $\chi^*$ in $\rho$, and hence the
multiplicity of $\lambda_{\chi}$ from \eqref{eigenvalue} as an
eigenvalue of $M$, is given by
\begin{equation}\label{eigenmult}
\frac{1}{|G_J|}\sum_{g\in G_J}\chi(g\inv)\sum_{J'\leq J, J'\in
\mathscr U(S)}|\mathrm{Fix}_L(e_{J'}ge_{J'})|\mu(J',J)
\end{equation}
where $\mu$ is the M\"obius function of $\mathscr U(S)$.

For the case where $S$ is an abelian group,~\eqref{eigenvalue} can
be found in the work of Diaconis~\cite{Diaconis}.  In this
situation, $L=S$ and the multiplicities are all one.  On the other
hand, if the semigroup is idempotent, then~\eqref{eigenvalue}
and~\eqref{eigenmult} reduce to the results of
Brown~\cite{Brown1,Brown2}.  If the maximal subgroups are trivial,
one obtains the results of~\cite{mobius1}.  Compare also
with~\eqref{DAvers}.

\section*{Acknowledgments}
This paper greatly benefitted from several e-mail conversations with
Mohan Putcha.  We also thank the anonymous referee for his many useful comments.  This paper would not have been possible if it had not been for Ariane Masuda, who first stirred my interest in the M\"obius function.

\bibliographystyle{amsplain}

\begin{thebibliography}{99}
\bibitem{Hopf}
M. Aguiar and R. Orellana, \textit{The Hopf algebra of uniform block
permutations}, J. Algebraic Combinatorics, to appear.

\bibitem{Aguiar}
M.~Aguiar and S.~Mahajan.
\newblock ``Coxeter Groups and {H}opf Algebras'', volume~23 of {\em Fields
  Institute Monographs}.
\newblock American Mathematical Society, Providence, RI, 2006.


\bibitem{Almeidabook}
J. Almeida, ``Finite Semigroups and Universal Algebra", World
Scientific, Singapore, 1995.


\bibitem{myrad}
J. Almeida, S.~W. Margolis, B. Steinberg and M.~V. Volkov,
\textit{Representation theory of finite semigroups, semigroup radicals
  and formal language
theory}, preprint 2004.


\bibitem{synch}
F. Arnold and B. Steinberg, \textit{Synchronizing groups and
automata}, Theor. Comp. Sci. \textbf{359} (2006), 101--110.

\bibitem{AshInv}
C.~J. Ash, \textit{Finite semigroups with commuting idempotents} J.
Austral. Math. Soc. Ser. A \textbf{43} (1987), 81--90.

\bibitem{ash}
C.~J. Ash,  \textit{Inevitable graphs: A proof of the type II
conjecture and some related decision procedures}, Internat. J.
Algebra Comput. \textbf{1} (1991), 127--146.
\bibitem{BHR}
P. Bidigare, P. Hanlon and D. Rockmore, \textit{A combinatorial
  description of the spectrum for the Tsetlin library and its
  generalization to hyperplane arrangements},
Duke Math. J. \textbf{99} (1999), 135--174.

\bibitem{Brown1}
K. Brown, \textit{Semigroups, rings, and Markov chains},
J. Theoret. Probab. \textbf{13} (2000), 871--938.

\bibitem{Brown2}
K. Brown, \textit{Semigroup and ring theoretical methods in
  probability}, in: ``Representations of Finite Dimensional Algebras
and Related Topics in Lie Theory and Geometry'', 3--26, (V. Dlab, ed.)
Fields Inst. Commun., \textbf{40}, Amer. Math. Soc., Providence, RI, 2004.

\bibitem{CP}
A.~H. Clifford and G.~B. Preston, ``The Algebraic Theory of
Semigroups'', Mathematical Surveys No. 7, AMS, Providence, RI, Vol. 1,
1961.

\bibitem{CR}
C.~W. Curtis and I.~Reiner, ``Representation theory of finite groups and associative algebras'',
Reprint of the 1962 original. AMS Chelsea Publishing, Providence, RI, 2006. xiv+689 pp.

\bibitem{Diaconis}
P. Diaconis, ``Group representations in probability and statistics",
Institute of Mathematical Statistics Lecture Notes---Monograph
Series, 11. Institute of Mathematical Statistics, Hayward, CA, 1988.

\bibitem{Eilenberg}
S. Eilenberg, ``Automata, Languages and Machines'',  Academic Press,
New York, Vol~A, 1974; Vol~B, 1976.



\bibitem{Green}
J.~A. Green, \emph{On the structure of semigroups}, Annals Math.
\textbf{54}, (1951), 163--172.

\bibitem{Hall}
M. Hall, ``Combinatorial theory'',
Reprint of the 1986 second edition. Wiley Classics Library. A
Wiley-Interscience Publication.
John Wiley \& Sons, Inc., New York, 1998.


\bibitem{HMPR}
K. Henckell, S.~W. Margolis, J.-E. Pin, and J. Rhodes, \textit{Ash's
type {II}
  theorem, profinite topology and Mal'cev products, Part {I}},
Internat. J. Algebra Comput. \textbf{1} (1991), 411--436.


\bibitem{Arbib} K. Krohn, J. Rhodes and B. Tilson,
\emph{Lectures on the algebraic theory of finite semigroups and
 finite-state machines}, Chapters 1, 5--9 (Chapter 6 with M.~A.
Arbib) of ``Algebraic Theory of Machines, Languages, and
Semigroups'', ed.\ M.~A.~Arbib, Academic Press, New York, 1968.

\bibitem{Lallement}
G. Lallement and M. Petrich, \emph{Irreducible matrix
representations of finite semigroups} Trans. Amer. Math. Soc.
\textbf{139} (1969), 393--412.


\bibitem{Lawson}
M.~V. Lawson, ``Inverse Semigroups: The theory of partial
symmetries'', World Scientific, Singapore,
1998.

\bibitem{MacLane}
S. MacLane, Categories for the Working Mathematician,
Springer-Verlag,
  New York, 1971.

\bibitem{McAlister}
D. B. McAlister,
\textit{Characters of finite semigroups},
J. Algebra \textbf{22} (1972), 183--200.

\bibitem{McAlistersurvey}
D.~B. McAlister, \emph{Representations of semigroups by linear
transformations I, II}, Semigroup Forum \textbf{2} (1971), 189--263;
ibid. \textbf{2} (1971), 283--320.

\bibitem{Munnalg}
W.~D. Munn,  \textit{On semigroup algebras},
Proc. Cambridge Philos. Soc. \textbf{51}, (1955). 1--15.

\bibitem{Munn}
W.~D. Munn,  \textit{Matrix representations of semigroups}, Proc.
Cambridge Philos. Soc. \textbf{53} (1957), 5--12.

\bibitem{Munnchar}
W.~D. Munn, \textit{The characters of the symmetric inverse semigroup},
Proc. Cambridge Philos. Soc. \textbf{53} (1957), 13--18.


\bibitem{Okninski}
J. Okni\'nski, ``Semigroup Algebras'', Monographs and Textbooks in Pure and Applied Mathematics, vol.~138, Marcel Dekker, Inc., New York, 1991.


\bibitem{Paterson}
A.~L.~T. Paterson, ``Groupoids, Inverse Semigroups, and their Operator
Algebras'', Birkh\"auser, Boston, 1999.

\bibitem{Poni}
I.~S. Ponizovski{\u\i}, \emph{On matrix representations of associative systems}, Mat. Sb. N.S.
\textbf{38}, (1956), 241--260.

\bibitem{Putcha}
M.~S. Putcha, ``Linear Algebraic Monoids'', London Mathematical
Society Lecture Note Series, \textbf{133}. Cambridge University Press,
Cambridge, 1988.

\bibitem{Putchachar}
M.~S. Putcha, \emph{Complex representations of finite monoids},
Proc.\ London Math.\ Soc. \textbf{73}  (1996),  623--641.

\bibitem{PutchaHecke}
M.~S. Putcha, \emph{Monoid Hecke algebras}, Trans. Amer. Math. Soc.
\textbf{349} (1997), 3517--3534.

\bibitem{quiver}
M.~S. Putcha, \emph{Complex representations of finite monoids. II.
Highest weight categories and quivers}, J. Algebra \textbf{205}
(1998), 53--76.

\bibitem{Putchaweights}
M.~S. Putcha, \emph{Semigroups and weights for group
representations},  Proc.\ Amer.\ Math.\ Soc.  \textbf{128} (2000),
2835--2842.

\bibitem{Renner}
L.~E. Renner, \textit{Analogue of the Bruhat decomposition for
  algebraic monoids},  J. Algebra  \textbf{101}  (1986),  303--338.

\bibitem{Rennerbook}
L.~E. Renner, ``Linear algebraic monoids", Encyclopaedia of
Mathematical Sciences, \textbf{134}, Invariant Theory and Algebraic
Transformation Groups, V. Springer-Verlag, Berlin, 2005.

\bibitem{rhodeschar}
J. Rhodes, \emph{Characters and complexity of finite semigroups}, J.
Combinatorial Theory  \textbf{6} (1969), 67--85.

\bibitem{qtheor}
J. Rhodes and B. Steinberg, \emph{The $\mathfrak q$-theory of finite
  semigroups}, Springer, to appear.

\bibitem{RhodesZalc}
J. Rhodes and Y. Zalcstein, \textit{Elementary representation and
character theory of finite semigroups and its application} in:
"Monoids and Semigroups with Applications" (J. Rhodes, ed.),
Berkeley, CA, 1989, 334--367, World Sci. Publishing, River Edge,
NJ, 1991.


\bibitem{RZ}
L. Ribes and P.~A. Zalesski{\u\i}, \textit{On the profinite topology
on a free
  group}, Bull. London Math. Soc. \textbf{25} (1993), 37--43.

\bibitem{Rosenblatt}
M. Rosenblatt, \textit{Stationary measures for random walks on
semigroups} in: ``Semigroups (Proc. Sympos., Wayne State Univ.,
Detroit, Mich., 1968)'', (K. Folley, ed.), pp. 209--220 Academic
Press, New York, 1969

\bibitem{sagan}
B.~E. Sagan, ``The symmetric group. Representations, combinatorial
algorithms, and symmetric functions," Second edition, Graduate Texts
in Mathematics, \textbf{203}, Springer-Verlag, New York, 2001.

\bibitem{Schutzrep}
M.-P. Sch\"utzenberger, \textit{$\overline{\D}$ repr\'esentation des
  demi-groupes},
C. R. Acad. Sci. Paris \textbf{244} (1957), 1994--1996.

\bibitem{Schutzambi}
M. P. Sch\"utzenberger, \textit{Sur le produit de concatenation
non ambigu},
 Semigroup Forum \textbf{13} (1976), 47--75.


\bibitem{Serre} J.-P. Serre,
``Linear representations of finite groups.'' Translated by Leonard
L. Scott. Graduate Texts in Mathematics, Vol. \textbf{42}.
Springer-Verlag, New York-Heidelberg, 1977.

\bibitem{Solomon1}
L. Solomon, \textit{The Burnside algebra of a finite group}, J.
Combinatorial Theory \textbf{2} (1967), 603--615.

\bibitem{Solomon3}
L.  Solomon, \textit{The Bruhat decomposition, Tits system and Iwahori ring
 for the monoid of matrices over a finite field},  Geom. Dedicata  \textbf{36}
 (1990),  15--49.

\bibitem{Solomon}
L. Solomon,
\emph{Representations of the rook monoid},
J. Algebra \textbf{256} (2002), 309--342.

\bibitem{Stanley}
R. Stanley,  ``Enumerative Combinatorics. Vol. 1'', Cambridge Studies
in Advanced Mathematics, vol.49, Cambridge University Press, 1997,
with a foreword by Gian-Carlo Rota, corrected reprint of the 1986 original.

\bibitem{mobius1}
B. Steinberg, \textit{M\"obius functions
  and semigroup representation theory}, J. Combinatorial Theory A \textbf{113} (2006), 866--881.

\bibitem{Therien}
H. Straubing and D. Th\'erien, \textit{Regular languages defined by
generalized first-order formulas with a bounded number of bound
variables}, Theory Comput. Syst. \textbf{36} (2003),  29--69.

\bibitem{zalc}
Y. Zalcstein, \emph{Studies in the representation theory of finite
semigroups},  Trans.\ Amer.\ Math.\ Soc. \textbf{161}  (1971),
71--87.
\end{thebibliography}

\end{document}